\documentclass[12pt,]{article}
\usepackage{pifont}
\usepackage{stmaryrd}
\renewcommand{\baselinestretch}{1.2}
\usepackage[normal]{subfigure}
\usepackage{graphicx}
\usepackage{amssymb}
\usepackage{amscd}
\usepackage{amsmath}
\usepackage{amsfonts}
\usepackage{amsbsy}
\usepackage{epsfig}

\begin{document}
\date{}

\newenvironment{ar}{\begin{array}}{\end{array}}
\baselineskip20pt \addtocounter{page}{0}
\title{\Large
Limit Cycle Bifurcations from Centers of Symmetric Hamiltonian Systems Perturbing by Cubic Polynomials
 \thanks{The project was supported by Key Disciplines of Shanghai Municipality (S30104), Slovenian Research Agency, and Slovene Human Resources Development and Scholarship Fund.} }

 \author{\small{{ZHAOPING HU $^{\mbox{a,b,}}$\thanks{Author for correspondence: zhaopinghu@shu.edu.cn.},}}
 \;BIN GAO $^{\mbox{a}}$,\;VALERY G. ROMANOVSKI $^{\mbox{b,c}}$ \\
 {\scriptsize $^{a}$ Department of
Mathematics,\ \ Shanghai University, Shanghai 200444, P. R. China}\\
{\scriptsize$^{b}$Center for Applied Mathematics and Theoretical Physics,
University of Maribor, SI-2000 Maribor, Slovenia}\\
 {\scriptsize$^{c}$ Faculty of Natural Science and Mathematics,\ \
University of Maribor, SI-2000 Maribor, Slovenia}}


\date{ }
\maketitle

\noindent{\scriptsize {{\bf Abstract.}
\small \baselineskip14pt \noindent In this paper, we consider some cubic near-Hamiltonian systems obtained from perturbing the symmetric cubic Hamiltonian system with two symmetric singular points by cubic polynomials. First, following Han [2012] we develop a method to study the analytical property of the Melnikov function near the origin for near-Hamiltonian system having the origin as its elementary center or nilpotent center. Based on the method, a computationally efficient algorithm is established to systematically compute the coefficients of Melnikov function. Then, we consider the symmetric singular points and present the conditions for one of them to be elementary center or nilpotent center. Under the condition for the singular point to be a center, we obtain the normal form of the Hamiltonian systems near the center. Moreover, perturbing the symmetric cubic Hamiltonian systems by cubic polynomials, we consider limit cycles bifurcating from the center using the algorithm to compute the coefficients of Melnikov function. Finally, perturbing the symmetric hamiltonian system by symmetric cubic polynomials, we consider the number of limit cycles near one of the symmetric centers of the symmetric near-Hamiltonian system, which is same to that of another center.}}\vspace{0.5cm}


\normalsize \baselineskip16pt

\section{Introduction}

There have been many studies on bifurcations of limit
cycles from elementary singular points, see Bautin [1952]-Yu \& Han [2004]. In general, there are
two types of such bifurcations leading to limit cycles: either by
perturbing a focus or by perturbing a center. Especially,
many authors study $C^\infty$ systems of the form
\begin{equation}\label{e11}
\dot{x}=H_y +\varepsilon p(x,y, \delta),\quad \dot{y}=-H_x
+\varepsilon q(x,y, \delta),
\end{equation}
where $H(x,y)$, $p(x,y,  \delta)$, $q(x,y, \delta)$ are $C^\infty$
functions, $\varepsilon\geq 0$ is small and $\delta\in D \subset
{\mathbf{ R}}^m$ is a vector parameter with $D$ compact.

When $\varepsilon=0$, system \eqref{e11} becomes
\begin{equation}\label{e12}
\dot{x}=H_y,  \quad  \dot{y}=-H_x,
\end{equation}
which is Hamiltonian system, and thus \eqref{e11} is called a
near-Hamiltonian system. We impose the following hypothesis

$(H)$ $H_x(x_0,y_0)=H_y(x_0,y_0)=0$\\
and denote $A(x_0,y_0)=\dfrac{\partial(H_y,-H_x)}{\partial(x,y)}(x_0,y_0)$ and $A=A(0,0)$.

If the hypothesis $(H)$ is satisfied, then system \eqref{e12} has a critical point at the point $(x_0,y_0)$. Suppose that the hypothesis $(H)$ is satisfied at the origin. Then it is said that the origin is $(i)$ a saddle if $det(A)<0$; $(ii)$ an elementary center if $det(A)>0$; $(iii)$ a nilpotent critical point if $det(A)=0$ and $A\neq0$.

\vspace{0.10in}

From now on, we always suppose that system \eqref{e12}  has a critical point at the origin. Without loss of generality, for $det(A)>0$ we may assume that the expansion of $H$ at the origin is of the form
\begin{equation}\label{ech}
H(x,y)=\dfrac{\omega}2(x^2+y^2)+\sum_{i+j\geq3}h_{ij}x^iy^j,\;\omega>0,
\end{equation}
for $det(A)=0$ and $A\neq0$ we may assume
\begin{equation}\label{nph}
H(x,y)=\dfrac{\omega}2\, y^2+\sum_{i+j\geq3}h_{ij}x^iy^j,\;\;\omega>0.
\end{equation}
We call \eqref{ech} or \eqref{nph} the normal form of Hamiltonian function of system \eqref{e12} having the origin as its elementary center or nilpotent critical point respectively. Similarly, we call the corresponding Hamiltonian system or near-Hamiltonian system the normal form of system \eqref{e12} or \eqref{e11} near the center.

For the function $H(x,y)$ of the normal form \eqref{ech} or \eqref{nph}, the implicit function theorem implies that there exists a unique analytic function $\varphi(x)=O(x^2)$ such that $H_y(x,\varphi(x))=0$ for $|x|$ small. Let
\begin{equation}
H^*_0(x)=H(x,\varphi(x))=\sum_{j\geq{k}}h_jx^j,\;h_k\neq0,\;k\geq2.
\end{equation}
 Recently, Han, et al. [2010] gave a complete classification of nilpotent critical points for Hamiltonian systems \eqref{e12} as follows.

\vspace{0.10in} \noindent {\bf Theorem 1.} \ \ {\it
The origin is $(i)$ a cusp if $k$ is odd; $(ii)$ a saddle if $k$ is even with $h_k<0$; and $(iii)$ a center if $k$ is even with $h_k>0$.}

\vspace{0.10in} \noindent {\bf Definition 1.} \ \ {\it
Let the Hamiltonian function $H(x,y)$ be of the form \eqref{ech} or \eqref{nph}. Then for system \eqref{e12} the origin is called a cusp of order $m$ if $k=2m+1$. It is called a nilpotent center of order
 $m$ (a nilpotent saddle of order $m$, respectively) if $k=2m+2$ and $h_k>0$ (if  $k=2m+2$ and $h_k<0$).}

 \vspace{0.10in}

Now suppose that the Hamiltonian system \eqref{e12} has an elementary center or a nilpotent center of order $p-1$ $(p\geq2)$ at the origin, namely the Hamiltonian function $H(x,y)$ has the form \eqref{ech} or \eqref{nph} and satisfies
\begin{equation}\label{hs0}
H^*_0(x)=\sum_{j\geq{2p}} h_j x^j,\;\;h_{2p}>0,\;\;p\geq1.
\end{equation}
Then, the Hamiltonian system \eqref{e12} has a family of periodic
orbits, given by
$$
L_h: H(x,y)=h , \quad   h\in (0,\beta)
$$
such that $L_h$ approaches the origin as $h\rightarrow 0$.

Take $h=h_0\in(0,\beta)$ and $A(h_0)\in L_{h_0}$. Let $l$ be a cross
section of system \eqref{e12} passing through $A(h_0)$. Then, for
$h$ near $h_0$ the periodic orbit $L_h$ has a unique intersection
point with $l$, denoted by $A(h)$, i.e., $A(h)=L_h\cap l$. Consider
the positive orbit $\gamma (h,\varepsilon,A)$ of system
\eqref{e11} starting from $A(h)$. Let $B(h,\varepsilon,A)$
denote the first intersection point of the orbit with $l$. Then we
have
\begin{equation}\label{e14}
H(B)-H(A)=\displaystyle\int_{AB}dH=\varepsilon[M(h,\delta)+O(\varepsilon)]
=\varepsilon \, F(h,\varepsilon,\delta),
\end{equation}
where
\begin{equation}\label{e15}
\begin{array}{ll}
M(h,\delta)&=\displaystyle\oint_{L_h}(H_yq+H_xp)dt  \vspace{0.1cm}\\
&=\displaystyle\oint_{L_h}(q dx-p dy)
=\displaystyle\int\!\!\!\!\!\!\int\limits_{H\leq h}(p_x+q_y) dx
  dy.
\end{array}\end{equation}
The functions $F(h,\varepsilon,\delta)$ and $M(h,\delta)$ in \eqref{e14} are called
bifurcation function and Melnikov function of system \eqref{e11} near the origin respectively. The resulting map from
$A(h)$ to $B(h, \varepsilon, \delta)$ is called a Poincare map of
system \eqref{e11}. Obviously, for small $\varepsilon$ system
\eqref{e11} has a limit cycle near the origin if and only if the
function $F(h,\varepsilon,\delta)$ has an isolated positive zero in $h$ near
$h=0$.

Based on the analytical property of the Melnikove function $M(h,\delta)$ and the number
of limit cycles near the origin by the function, we have the
following theorems (see Han, Jiang $\&$ Zhu [2008]).

\vspace{0.10in} \noindent {\bf Theorem 2.} \ \ {\it Let system \eqref{e12}, where $H(x,y)$ has the form \eqref{ech} or \eqref{nph}, satisfy \eqref{hs0}. Then, we have
\begin{equation}
M(h,\delta)=h^{\frac{p+1}{2p}} \sum_{j\geq0}b_j(\delta) h^{\frac{j}p}.
\end{equation}
In particular, when p=1 ( that is when system \eqref{e12} has an elementary center at the origin), we have
\begin{equation*}
M(h,\delta)=h \sum_{j\geq0}b_j(\delta) h^j.
\end{equation*}
}

\vspace{0.10in} \noindent {\bf Theorem 3.}\ \ {\it Under the
condition of Theorem $2$,  if there exist $k\geq 1,\,\,\delta_0\in
D$ such that $b_k(\delta_0)\neq0$ and
$$
b_j(\delta_0)=0,\,\,\,\,j=0,1,\cdots,k-1,
\hskip 6mm \hbox{\em{det}}\frac{\partial(b_0,\cdots,b_{k-1})}
{\partial(\delta_1,\cdots,\delta_k)}(\delta_0)\neq0,
$$
where $\delta=(\delta_1,\cdots,\delta_m),$ $m\geq k$, then there
exist a constant $\varepsilon_0 >0$ and a neighborhood $V$ of the
origin such that for some $0<|\varepsilon|<\varepsilon_0$ and
$|\delta-\delta_0|<\varepsilon_0$  \eqref{e11} has $k$ limit cycles in $V_1$. }

\vspace{0.10in} \noindent {\bf Theorem 4.} \ {\it Consider the
near-Hamiltonian system \eqref{e11}, where $H(x,y,a)$ with $a\in
{\mathbf{ R}}^n$ satisfies \eqref{hs0} and the functions $p$ and $q$
are linear in $\delta \in {\mathbf{ R}}^m$. Suppose there exist
integer $k>0$ and
$\delta_0=(\delta_{10},\cdots,\delta_{m0})\in{\mathbf{ R}}^m$, $a_0
\in{\mathbf{ R}}^n$ such that
\begin{equation}\label{e222}
b_j(\delta_0,a_0)=0, \ j=0,\cdots, k-1,\ \det
\frac{\partial(b_0,\cdots,b_{k-1})}
 {\partial(\delta_1,\cdots,\delta_k)}(a_0)\neq0,
 \end{equation}
and
\begin{equation}\label{e223}
b_{k+j}|_{(\delta_1,\cdots,\delta_k)=\xi(\delta_{k+1},\cdots,\delta_m,\,a)}
=L_j(\delta_{k+1},\cdots,\delta_m)\Delta_j(a), \ j=0,\cdots,n,
\end{equation}
where
$$(\delta_1,\cdots,\delta_k)=\xi(\delta_{k+1},\cdots,\delta_m,a)$$
is the unique solution to $b_j=0,\,j=0,\cdots,k-1$ for a point $a$ near $a_0$,
\begin{equation}\label{e224}
\begin{array}{c}L_j(\delta_{k+1,0},\cdots,\delta_{m0})\neq 0, \ j=0,\cdots,
n,\\ [1.0ex]
 \Delta_j(a_0)=0, \  j=0,\cdots, n-1,\  \Delta_n(a_0)\neq 0,
\end{array}
\end{equation}
and
\begin{equation}\label{e225}
\det \frac{\partial(\Delta_0,\cdots,\Delta_{n-1})}
 {\partial(a_1,\cdots,a_n)}(a_0)\neq0.
 \end{equation}
Then, for some $(\varepsilon,\delta,a)$ near $(0,\delta_0,a_0)$,
system \eqref{e11} has $k+n$ limit cycles near the origin.}

\vspace{0.10in}

In light of the above theorems, a key step in studying the small-amplitude limit cycle bifurcations of system \eqref{e11} is to find an efficient
method to compute the coefficients $b_l$. For p=1, the formulas for the first
three coefficients $b_j(\delta),\,j=0,1,2$ were obtained by Hou \&
Han [2006] by using the double integral in \eqref{e15}. Recently, Han, Yang $\&$ Yu [2009] and Han[2012] have developed a new approach to prove Theorem 2 and have established a computationally efficient algorithm to systematically compute $b_j(\delta),\,j=0,\,1,\,2,\,3,\,\cdots$ for $p=1$ and $p>1$ respectively. Following Han, Yang $\&$ Yu [2009] and Han[2012], we prove Theorem 2 and then also establish a computationally efficient algorithm to systematically compute $b_j(\delta),\,j=0,\,1,\,2,\,3,\,\cdots$ for arbitrary positive integer $p$.

In this paper, to illustrate the efficiency of the approach, we consider the following symmetric cubic Hamiltonian system\\
$(S)$: Hamiltonian system \eqref{e12}, with Hamiltonian function of the form
\begin{equation*}
H(x,y)=\sum_{i+j=2,4} h_{ij} x^i y^j,
\end{equation*}
taking $(0,\pm1)$ as its critical points,

We will consider the conditions for the symmetric points to be elementary centers or nilpotent centers. Under the condition for one of the symmetric points to be a center, perturbing system $(S)$ by cubic polynomials one can obtain a cubic near-Hamiltonian system. Then, by suitable change of the variables $(x,y)$ and $t$, we gave out the normal form of the above near-Hamiltonian system near the center $(0,1)$ (or $(0,-1)$). Obviously, the number of limit cycles of the near-Hamiltonian system near the center $(0,1)$ (or $(0,-1)$) is the same as that of the corresponding system of the normal form near the origin. As the result, we can study the number of limit cycles near the center $(0,1)$ (or $(0,-1)$) by studying the Melnikov of the normal form near the origin. Moreover, perturbing the symmetric hamiltonian system $(S)$ by symmetric polynomials, one can obtain a symmetric near-Hamiltonian system. Then, if one found some limit cycles near one of the symmetric centers, then the system has the same number of limit cycles near the other center. Hence, it is easier to find some limit cycles of the system near both centers than to find the same number of limit cycles near one center.

This paper is organized as follows. In the next section, we give a proof of Theorem 2, and list an efficient algorithm for any positive integer $p$ based on this proof as the appendix. In section 3, we study the critical points $(0,\pm1)$ of system $(S)$ and list all conditions for one of the singular points to be an elementary center or a nilpotent center. Moreover, near the center we present the corresponding normal form of the Hamiltonian system $(S)$. Then, in section 4, perturbing system $(S)$ by cubic polynomials, we obtain the corresponding near-Hamiltonian system. Applying the theorems and the program to the normal form of the above near-Hamiltonian system, we estimate the number of limit cycles bifurcating from the center $(0,1)$ or $(0,-1)$. In section 5, we discuss the number of limit cycles near both the symmetric centers of the symmetric cubic near-Hamiltonian systems obtained from perturbing system $(S)$ by symmetric cubic polynomials.

\section{ Proof of Theorem 2}

Recently, Han, Yang $\&$ Yu [2009] developed a new approach to prove Theorem 2 for $p=1$ and then established an algorithm based on its proof. The algorithm was implemented in the computer algebra system Maple. Earlier, Han, Jiang $\&$ Zhu [2008] proved Theorem 2 for $p>1$. Then following Han, Yang $\&$ Yu [2009], Han [2012] gave a new proof. In this paper, for convenience and as a preliminary, we repeat the proof of Theorem 2 for all positive integer $p\geq1$ following Han [2012], and then establish an algorithm based on the proof and provide an implementation in computer algebra system Mathematica.

To prove Theorem 2, we first introduce a change of
variables to make the form of the Hamiltonian function simpler. By
\eqref{ech} and \eqref{nph} using the implicit function theorem, one can show
that there exists a unique $C^{\infty}$ function $\varphi(x)$ such that
$H_y(x,\varphi(x))=0$ for $|x|$ small. Thus, we can write
\begin{equation}\label{e21}
\varphi(x)=\sum_{j\geq 2}e_jx^j.
\end{equation}

By introducing a new variable $v=y-\varphi(x)$, system \eqref{e11}
can be rewritten as
\begin{equation}\label{e22}
\begin{array}{l}
\dot{x}=H_v^{*}(x,v)+\varepsilon
p^*(x,v, \delta), \vspace{0.2cm}\\
\dot{v}=-H_x^{*}(x,v)+\varepsilon q^*(x,v, \delta),
\end{array}\end{equation}
where
\begin{equation}\label{e23}
\begin{array}{ll}
H^*(x,v)&=H(x,v+\varphi(x)),\\[1.0ex]
p^*(x,v, \delta)&= p(x,v+\varphi(x), \delta),\\[1.0ex]
q^*(x,v, \delta)&= q(x,v+\varphi(x), \delta)-\varphi^{'}(x)p^*(x,v,
\delta).
\end{array}
\end{equation}
Noting that $H_y(x,\varphi(x))=0$, we have
\begin{equation}\label{hxv}
H^*(x,v)=H^*_0(x)+\sum_{j\geq1}H^*_j(x)v^{j+1}=H^*_0(x)+v^2 \tilde{H}(x,v),
\end{equation}
where $\tilde{H}(0,0)=\dfrac{\omega}2>0$ and
\begin{equation}\label{hst}
\begin{array}{ll}
H^*_0(x)&=H(x,\varphi(x))=\displaystyle\sum_{j\geq2}h_j x^j,  \vspace{0.1cm}\\[1.0ex]
H^*_j(x)&=\dfrac{1}{(j+1)!} \dfrac{\partial^{j+1}H}{\partial{y}^{j+1}} (x,\varphi(x)).
\end{array}
\end{equation}

Let \eqref{hs0} be satisfied. Then it follows from \eqref{hxv} that there exist a family of periodic
orbits surrounding the origin defined by the equation $H(x,y)=h$ or
$H^*(x,v)=h$ for $h>0$ small. Let $w=\sqrt{h-H^*_0(x)}$ and suppose that the region surrounding by the closed curve $H^*(x,v)=h$ can be expressed as the form $\{(x,v)|x_2(h)\leq{x}\leq{x_1(h)},\;v_2(x,w)\leq{v}\leq{v_1(x,w)}\}$. Then, from Han, Yang $\&$ Yu [2009] or Han [2012], we have

\vspace{0.10in} \noindent {\bf Lemma 1.}\ \ {\it
$(i)$ The equation $H^*(x,v)=h$ has exactly two $C^{\infty}$ solutions $v_1(x,w)>0$ and
$v_2(x,w)<0$  in $v$ satisfying
$$
v_1(x,w)=\sqrt{2}w(1+O(|x,w|)),\quad v_2(x,w)=v_1(x,-w).
$$
$(ii)$ $x_1(h)>0$ and $x_2(h)<0$ are the two solutions to equation $H^*_0(x)=h$.}

\vspace{0.10in}

Then, it is obvious that
\begin{equation}\label{e26}
\begin{array}{ccl}
M(h,\delta)\!\!\!&=&\!\!\!\displaystyle\oint_{H^*(x,v)=h}q^*dx-p^*dv \vspace{0.2cm}\\
&=&\!\!\!\displaystyle\int\!\!\!\!\int_{H^*\leq h}(p^*_x+q^*_v)dxdv\vspace{0.2cm}\\
&=&\!\!\!\displaystyle\int_{x_2(h)}^{x_1(h)}\displaystyle\int_{v_2(x,w)}^{v_1(x,w)}(p^*_x+q^*_v)dvdx.
\end{array}\end{equation}

Let
\begin{equation}\label{e27}
\bar{q}(x,v,\delta)=\displaystyle\int_{0}^{v}(p^*_x+q^*_v)dv=q^*(x,v,\delta)-q^*(x,0,\delta)+\int_0^vp^*_x(x,u,\delta)du.
\end{equation}
Then, we have
\begin{equation}\label{e29}
\bar{q}(x,v,\delta)=v\sum_{i+j\geq 0}\bar{b}_{ij}x^iv^j=\sum_{j\geq
1}q_j(x)v^j,
\end{equation}
 where
\begin{equation}\label{e210}
\begin{array}{rcl}
q_{j+1}(x)&=&\displaystyle\frac{1}{(j+1)!}\frac{\partial^j}{\partial
v^j}(p_x^*+q_v^*)|_{\varepsilon=v=0}\vspace{0.2cm}\\
&=&\displaystyle\frac{1}{(j+1)!}\frac{\partial^j}{\partial
y^j}(p_x+q_y)(x,\varphi(x),0,\delta)\vspace{0.2cm}\\
&=&\displaystyle\sum_{i\geq 0}\bar{b}_{ij}x^i,\quad j\geq 0.
\end{array}\end{equation}

Further, we write
\begin{equation}\label{e212}
v_1(x,w)=\sum_{j\geq 1}a_j(x)w^j.
\end{equation}
By Lemma 1 we have
\begin{equation*}
\begin{array}{rcl}
w^2&=&h-H_0^*(x)=H^*(x,v)-H_0^*(x)=(v_1(x,w))^2\displaystyle\sum_{j\geq
1}H_j^*(x)[v_1(x,w)]^{j-1}\vspace{0.2cm}\\
&=&a_1^2(x)H_1^*(x)w^2+(2a_1(x)a_2(x)H_1^*(x)+a_1^3(x)H_2^*(x))w^3+\vspace{0.2cm}\\
&&[(a_2^2(x)+2a_1(x)a_3(x))H_1^*(x)+3a_1^2(x)a_2(x)H_2^*(x)+a_1^4(x)H_3^*(x)]w^4+\cdots.
\end{array}
\end{equation*}
Equaling the coefficients of $\, w^j$ in the above identity, we have
\begin{equation}\label{e213}
a_1(x)=\frac{1}{\sqrt{H_1^*(x)}},\
a_2(x)=-\frac{H_2^*(x)}{2(H_1^*(x))^2},\
a_3(x)=-\frac{1}{8(H_1^*)^{\frac{7}{2}}}[4H_1^*H_3^*-5(H_2^*)^2],\
\cdots. \end{equation}

Then, it follows from \eqref{e26} that
$$
M(h,\delta)=\int_{x_2(h)}^{x_1(h)}[\bar{q}(x,v_1(x,w))-\bar{q}(x,v_2(x,w))]dx.
$$
By Lemma 1, the function $\bar{q}(x,v_1)-\bar{q}(x,v_2)$ is odd in
$w$. Thus, we can write
\begin{equation}\label{e214}
\bar{q}(x,v_1)-\bar{q}(x,v_2)=\sum_{j\geq 0}\bar{q}_j(x)w^{2j+1},
\end{equation}
and hence,
$$
M(h,\delta)=\sum_{j\geq
0}\int_{x_2(h)}^{x_1(h)}\bar{q}_j(x)w^{2j+1}dx.
$$

In order to compute the above integral, we change the limits of
integration. Let $\psi(x)={\rm sgn}(x)[H_0^*(x)]^{\frac{1}{2p}}$.
Then, by \eqref{hs0} the function $\psi$ is $C^{\infty}$ for small
$|x|$ with $\psi'(0)=h_{2p}^{\frac{1}{2p}}>0$. Therefore, we may introduce
the new variable $u=\psi(x)$ to obtain
$$
\begin{array}{rcl}
M(h,\delta)&=&\displaystyle\sum_{j\geq
0}\int_{-h^{\frac{1}{2p}}}^{h^{\frac{1}{2p}}}\tilde{q}_j(u)w^{2j+1}du\vspace{0.1cm}\\
&=&\displaystyle\sum_{j\geq
0}\int_{0}^{h^{\frac{1}{2p}}}[\tilde{q}_j(u)+\tilde{q}_j(-u)]w^{2j+1}du,
\end{array}
$$
where $w=\sqrt{h-u^{2p}}$ and
\begin{equation}\label{e215}
\tilde{q}_j(u)=\left.\frac{\bar{q}_j(x)}{\psi^{'}(x)}\right|_{x=\psi^{-1}(u)}.
\end{equation}
It is easy to see that we can assume
\begin{equation}\label{e216}
\tilde{q}_j(u)+\tilde{q}_j(-u)=\sum_{i\geq 0}r_{ij}u^{2i}.
\end{equation}
Then,
\begin{equation}\label{e217}
M(h,\delta)=\sum_{i+j\geq 0}r_{ij}I_{ij}(h),
\end{equation}
where
$$
I_{ij}(h)=\displaystyle\int_0^{h^{\frac{1}{2p}}}u^{2i}w^{2j+1}du=\displaystyle\int_0^{h^{\frac{1}{2p}}}u^{2i}(h-u^{2p})^j\sqrt{h-u^{2p}}du.
$$

\vspace{0.10in}
\noindent {\bf Lemma 2.}\ \ {\it Let
\begin{equation}\label{e218}
\beta_{ij}=\int_0^1v^{\frac{2i}p}(1-v^2)^j\sqrt{1-v^2}dv.
\end{equation}
Then
$$
I_{ij}(h)=\beta_{ij}h^{\frac{p+1}{2p}+\frac{i+pj}p},\quad h>0, \,\,\, 0<\beta_{ij}<1.
$$}

\vspace{-0.15in} \noindent
\emph{Proof.}\ \ Introducing
$v=u^p/h^{\frac{1}{2}}$, we have
$$
I_{ij}(h)=h^{\frac{p+1}{2p}+\frac{i+pj}p}\int_0^1v^{\frac{2i}p}(1-v^2)^j\sqrt{1-v^2}dv.
$$
This ends the proof.
\vspace{0.10in}

Now by \eqref{e217} and Lemma 2, we have
\begin{equation}\label{e219}
M(h,\delta)=\displaystyle h^{\frac{p+1}{2p}}\sum_{i+j\geq 0}r_{ij}\beta_{ij}h^{\frac{i+pj}p}
=\displaystyle h^{\frac{p+1}{2p}}\sum_{l\geq 0}b_l(\delta)h^{l},
\end{equation}
where
\begin{equation}\label{e220}
b_l(\delta)=\sum_{i+pj=l}r_{ij}\beta_{ij}.
\end{equation}

Finally based on \eqref{e214}, \eqref{e215} and \eqref{e216}, we
see that if system \eqref{e11} is analytic, then the series
$\sum\limits_{i,j\geq 0}r_{ij}u^{2i}w^{2j+1}$ is convergent for
$(u,w)$ near the origin. Then it follows that the series
$\sum\limits_{i,j\geq 0}|r_{ij}|\mu^{i+j}$ is convergent for some
constant $\mu>0$, and hence  by \eqref{e219} $M(h,\delta)$ is
analytic in $h$.

This completes the proof of Theorem 2.

\section{The symmetric critical points of symmetric Hamiltonian system $(S)$}

In this section, we study the symmetric critical points $(0,\pm1)$ of symmetric cubic Hamiltonian system $(S)$. Obviously, we have
$H_y(0,\pm1)=H_x(0,\pm1)=0$, i.e.,
$$\pm(2h_{02}+4h_{04})=\pm(h_{11}+h_{13})=0.$$
Thus, we can take
\begin{equation}\label{csp}
h_{04}=-h_{02}/2,\;h_{13}=-h_{11}.
\end{equation}

For the symmetry, we need only to consider the critical point $(0,1)$. By introducing a new variable $v=y-1$ and then taking $v$ as $y$, we have
\begin{equation}\label{bhs}
\dot{x}=\bar{H}_y(x,y),\;\;\dot{y}=-\bar{H}_x(x,y),
\end{equation}
where $\bar{H}(x,y)=H(x,y+1)$. Truncating the constant term of the function $\bar{H}$, we have
\begin{equation}\label{bhf}
\begin{array}{ll}
\bar{H}(x,y)&=(h_{20}+h_{22})x^2-2h_{11}xy-2h_{02}y^2+h_{31}x^3+2h_{22}x^2y-3h_{11}xy^2\vspace{0.2cm}\\
&-2h_{02}y^3+h_{40}x^4+h_{31}x^3y+h_{22}x^2y^2-h_{11}xy^3-\dfrac12h_{02}y^4.
\end{array}
\end{equation}

As a result, the point $(0,1)$ of system $(S)$ corresponds to the origin of system \eqref{bhs}. If \eqref{csp} is satisfied, then system \eqref{bhs} has a critical point at the origin. Moreover, for system \eqref{bhs} we have
\begin{equation}
A=\left(\begin{array}{cc}
         -2h_{11} & -4h_{02} \vspace{0.2cm}\\
         -2(h_{20}+h_{22}) & 2h_{11}
       \end{array}\right)
\end{equation}
and
\begin{equation}
det(A)=-4h^2_{11}-8h_{02}(h_{20}+h_{22}).
\end{equation}

Then, the origin is $(i)$ a saddle if $det(A)<0$; $(ii)$ an elementary center if $det(A)>0$; $(iii)$ a nilpotent critical point if $det(A)=0$ and $A\neq0$.
\vspace{0.2cm}

If $\Delta=det(A)>0$, which means that $h_{02}(h_{20}+h_{22})\neq0$, then system \eqref{bhs} has an elementary center at the origin. Let $\omega=\sqrt{\Delta}>0$. Without loss of generality, we can let $h_{20}+h_{22}=1/2$ and $\omega=1$ (otherwise, we can introduce a suitable rescaling of $(x,y)$ and time $t$). Therefore, we obtain
$$h_{20}=-h_{22}+\dfrac12\,,\;\,h_{02}=-h^2_{11}-\dfrac14\,.$$
Then, introducing the change of variables
\begin{equation}
u=-y,\;v=x-2 h_{11}y,
\end{equation}
and then taking $(u,v)$ as $(x,y)$, system \eqref{bhs} becomes
\begin{equation}\label{3hs1}
\dot{x}=\dfrac{\partial{H_1}}{\partial{y}}(x,y),\; \dot{y}=-\dfrac{\partial{H_1}}{\partial{x}}(x,y),
\end{equation}
where
\begin{equation}\label{hf1}
H_1(x,y)=\dfrac12(x^2+y^2)+\sum_{3\leq{i+j}\leq4}\bar{h}_{ij}x^i y^j,
\end{equation}
and
\begin{equation*}
\begin{array}{ll}
\bar{h}_{03}&= h_{31},\;\bar{h}_{12}= -2 h_{22} - 6 h_{11} h_{31},\;
\bar{h}_{21}= -3 h_{11} +8 h_{11} h_{22} +12 h^2_{11} h_{31},\vspace{0.2cm}\\
\bar{h}_{30}&= -\dfrac12 - 4(2 h_{22} +2 h_{11} h_{31}-1) h^2_{11},\;
\bar{h}_{04}=h_{40},\;
\bar{h}_{13}= -h_{31} -8 h_{11} h_{40},\vspace{0.2cm}\\
\bar{h}_{22}&= h_{22} + 6 h_{11} h_{31}+24 h^2_{11} h_{40},\;
\bar{h}_{31}= h_{11}-4 h_{11} h_{22} -12 h^2_{11} h_{31} -32 h^3_{11} h_{40},\vspace{0.2cm}\\
\bar{h}_{40}&=\dfrac18-\dfrac32h^2_{11} + 4 h^2_{11} h_{22} + 8 h^3_{11} h_{31} + 16 h^4_{11} h_{40}.
\end{array}
\end{equation*}

Obviously, if the symmetric Hamiltonian system $(S)$ has an elementary center at $(0,1)$, then the normal form of system $(S)$ near $(0,1)$ should be system \eqref{3hs1}. For the symmetry, system $(S)$ has an elementary center at $(0,-1)$ and is of the same normal form \eqref{3hs1} at the same time.
\vspace{0.2cm}

Now, suppose that system \eqref{bhs} has a nilpotent critical point at the origin. In other words, we have
\begin{equation}\label{cncp}
-4h^2_{11}-8h_{02}(h_{20}+h_{22})=0,\;h^2_{11}+h^2_{02}+(h_{20}+h_{22})^2\neq0.
\end{equation}

If $h_{02}=0$, then $h_{11}=0$ and $h_{20}+h_{22}\neq0$. Thus, system \eqref{bhs} becomes
\begin{equation*}
\begin{array}{ll}
\dot{x}&=x^2(2h_{22}-h_{31}x-2 h_{22}y)\vspace{0.3cm}\\
\dot{y}&=-x[2(h_{20}+h_{22})+3 h_{31}x+4 h_{22}y+4 h_{40}x^2+3h_{31}x y+2h_{22}y^2],
\end{array}
\end{equation*}
which means that all the points on the line $x=0$ are singular. Then the origin is not an isolated singular point.

Therefore, we can assume that $h_{02}\neq0$. Without loss of generality, suppose $h_{02}=\dfrac12$ (otherwise, we can introduce suitable rescaling of $(x,y)$). Then, from \eqref{cncp} we can take
\begin{equation}
h_{22}=-h^2_{11}-h_{20}.
\end{equation}
Introducing new variables as
\begin{equation}
u=-\dfrac12 x,\;\,v=h_{11} x+y,
\end{equation}
and then taking $(u,v)$ as $(x,y)$, system \eqref{bhs} becomes
\begin{equation}\label{bhs2}
\dot{x}=\dfrac{\partial{\bar{H}_2}}{\partial{y}}(x,y),\;\;\dot{y}=-\dfrac{\partial{\bar{H}_2}}{\partial{x}}(x,y),
\end{equation}
where
\begin{equation}\label{bhf2}
\bar{H}_2(x,y)=\frac{1}{2} y^2+\sum_{3\leq{i+j}\leq4}\tilde{h}_{ij}x^i y^j,
\end{equation}
and
\begin{equation}
\begin{array}{ll}
\tilde{h}_{30}&=8h_{11}h_{20}+4h_{31},\;\;\tilde{h}_{21}=-2h^2_{11}+4h_{20},\;\;\tilde{h}_{12}=\tilde{h}_{13}=0,\vspace{0.2cm}\\
\tilde{h}_{03}&=\dfrac{1}2,\;\;\tilde{h}_{40}=2h^2_{11}(h^2_{11}+4h_{20})+8(h_{11}h_{31}-h_{40}),\vspace{0.2cm}\\
\tilde{h}_{31}&=8 h_{11}h_{20}+4 h_{31},\;\;\tilde{h}_{22}=-h^2_{11}+2 h_{20},\;\tilde{h}_{04}=\dfrac{1}8\,.\vspace{0.2cm}
\end{array}
\end{equation}

Using the mathematica code we compute
\begin{equation}
\begin{array}{ll}
e_2&=-4h_{20}+2 h^2_{11}\;,\vspace{0.2cm}\\
e_3&=-4(2h_{11}h_{20}+ h_{31})\;,\vspace{0.2cm}\\
e_4&=-2 (h^2_{11}-2 h_{20})^2\;,\vspace{0.2cm}\\
e_5&=16 (h^2_{11}-2 h_{20})(2 h_{11}h_{20}+ h_{31})\;,\vspace{0.2cm}\\
e_6&=4 (h^6_{11}-6 h^4_{11} h_{20}-12 h^2_{11} h^2_{20}-8h^3_{20}-24 h_{11}h_{20}h_{31}-6h^3_{31})\;,\vspace{0.2cm}\\
&\;\;\;\;\;\;\;\cdots\cdots
\end{array}
\end{equation}
and
\begin{equation}
\begin{array}{ll}
h_3&=4 (2 h_{11}h_{20}+ h_{31})\;,\vspace{0.2cm}\\
h_4&=8 h_{11}(2 h_{11}h_{20}+ h_{31})-8 (h^2_{20}+ h_{40})\;,\vspace{0.2cm}\\
h_5&=8 (h^2_{11}-2 h_{02})(2 h_{11}h_{20}+ h_{31})\;,\vspace{0.2cm}\\
h_6&=-8 (2 h_{11}h_{20}+h_{31})^2\;,\vspace{0.2cm}\\
&\;\;\;\;\;\;\;\cdots\cdots
\end{array}
\end{equation}

Introduce the following conditions:

$(A)$ $h_{31}\neq -2 h_{11}h_{20}$;

$(B)$ $h_{31}=-2 h_{11}h_{20}$, $h^2_{20}+ h_{40}>0$;

$(C)$ $h_{31}=-2 h_{11}h_{20}$, $h^2_{20}+ h_{40}<0$;

$(D)$ $h_{31}=-2 h_{11}h_{20}$, $h_{40}=-h^2_{20}$.

\vspace{0.10in}
Then, we can prove the following statement.

\vspace{0.10in} \noindent {\bf Theorem 5.} \ {\it For system \eqref{bhs2}, the origin is an isolated nilpotent singular point if and only if one of conditions $(A),\;(B),\;(C)$ listed above holds. Further, the origin is

 $(i)$ a cusp of order 1 if and only if $(A)$ holds,

 $(ii)$ a nilpotent saddle of order 1 if and only if $(B)$ holds,

 $(iii)$ a nilpotent center of order 1 if and only if $(C)$ holds.}

\vspace{0.10in} \noindent \emph{Proof.} By Definition 1, the origin is a cusp of order 1 under $(A)$. If $(B)$ or $(C)$ holds, then $h_3=0,\;h_4<0$ or $h_3=0,\;h_4>0$, and the origin is a nilpotent saddle or center of order 1. If $(D)$ holds, then we have
$$\bar{H}_2(x,y)=\dfrac{1}{8}(y^2+2 y-4h^2_{11}x^2+8 h_{20}x^2)^2,$$
which means that all points on the curve
$$y^2+2 y-4h^2_{11}x^2+8 h_{20}x^2=0$$
are singular. This ends the proof.

\vspace{0.10in}

Assume that the origin is a nilpotent center of order 1 of system \eqref{bhs2}, i.e., $h_{31}=-2 h_{11}h_{20}$ and $h^2_{20}+h_{40}<0$. Then system \eqref{bhs2} becomes
\begin{equation}\label{3hs2}
\dot{x}=\dfrac{\partial{H_2}}{\partial{y}}(x,y),\; \dot{y}=-\dfrac{\partial{H_2}}{\partial{x}}(x,y),
\end{equation}
where
\begin{equation}\label{hf2}
H_2(x,y)=\dfrac12y^2+2Ax^2y+\dfrac12y^3+Bx^4+Ax^2y^2+\dfrac18y^4
\end{equation}
with
\begin{equation*}
A=2 h_{20}-h^2_{11},\;B=2h^4_{11}-8h^2_{11}h_{20}-8 h_{40}
\end{equation*}
satisfying $B>2\,A^2$.

Similarly, if the symmetric Hamiltonian system $(S)$ has a nilpotent center of order 1 at $(0,1)$, then the normal form of system $(S)$ near $(0,1)$ should be system \eqref{3hs2}. For the symmetry, system $(S)$ has a nilpotent center of order 1 at $(0,-1)$ and is of the same normal form \eqref{3hs2}.

\section{Bifurcation from the centers of systems $(S)$}

From section 3, we know that the symmetric Hamiltonian system $(S)$ can take $(0,\pm1)$ as elementary centers or nilpotent centers of order 1. Under the conditions for one of the symmetric critical points to be a center of system $(S)$, we perturb it by cubic polynomials to obtain a near-Hamiltonian system and then study limit cycles bifurcating from the center $(0,1)$ or $(0,-1)$.

Before starting to consider limit cycles bifurcating from a center of the near-Hamiltonian systems, we present some results on how to calculating the coefficients of the Melnikov function near the center. Suppose that system \eqref{e12} has a center at the point $(x_0,y_0)$ and $H(x_0,y_0)=h_0$. Then, there exist a family of periodic orbits given by
$$
L_h:\;H(x,y)=h,\;\;h\in(h_0,h_0+\beta),\;\beta>0.
$$
To transform the linear part of the Hamiltonian system at a center into a normal form, one often introduces suitable linear change of variables or rescaling of the time. This procedure may cause a change of the Melnikov function. More precisely, making a linear change of the form
\begin{equation*}
u=a(x-x_0)+b(y-y_0),\;\;v=c(x-x_0)+d(y-y_0),
\end{equation*}
and rescaling time $\tau=k\ t$, where $D=ad-bc\neq0$, system \eqref{e11} becomes
\begin{equation}\label{ths}
\dfrac{du}{d\tau}=\tilde{H}_v+\varepsilon\ \tilde{p},\;\,\dfrac{dv}{d\tau}=-\tilde{H}_u+\varepsilon\ \tilde{q},
\end{equation}
where
\begin{equation*}
\begin{array}{ll}
\tilde{H}(u,v)&=\dfrac{D}{k}[H(x,y)-h_0],\vspace{0.2cm}\\
\tilde{p}(u,v,\delta)&=\dfrac{1}{k}\ [a\ p(x,y,\delta)+b\ q(x,y,\delta)],\vspace{0.2cm}\\
\tilde{q}(u,v,\delta)&=\dfrac{1}{k}\ [c\ p(x,y,\delta)+d\ q(x,y,\delta)],\vspace{0.2cm}
\end{array}
\end{equation*}
$$
x=x_0+\dfrac1D(du-bv),\;y=y_0+\dfrac1D(-cu+av).
$$
Let $\tilde{M}$ denote the Melnikov function of system \eqref{ths} at the origin, which is given by
$$
\tilde{M}(h,\delta)=\oint_{\tilde{H}(u,v)=h}\tilde{q}\ du-\tilde{p}\ dv|_{\varepsilon=0},\;h\in(0,\beta').
$$
Then, it is easy to see that
$$
\begin{array}{ll}
\tilde{M}(h,\delta)&=sgn(k)\oint_{H(x,y)=h_0+kh/D}[(a\tilde{q}-c\tilde{p})dx+(b\tilde{q}-d\tilde{p})dy]|_{\,\varepsilon=0}\vspace{0.2cm}\\
&=sgn(k)\dfrac{D}{k}\oint_{H(x,y)=h_0+kh/D}[qdx-pdy]|_{\,\varepsilon=0}\vspace{0.2cm}\\
&=\dfrac{D}{|k|}M(h_0+\dfrac{k}{D}h,\delta),\;h\in(0,\beta'),\vspace{0.2cm}
\end{array}
$$
where $M(h,\delta)$ is the Melnikov function of system \eqref{e11} near the center $(x_0,y_0)$ of the corresponding Hamiltonian system \eqref{e12}.
Therefore, for $h_0<h<h_0+\beta$ we have
\begin{equation}
M(h,\delta)=\dfrac{|k|}{D}\tilde{M}(\dfrac{D}{k}(h-h_0,\delta)).
\end{equation}

By Theorem 2, we have

\vspace{0.10in} \noindent {\bf Theorem 6.} \ \ {\it Let system \eqref{ths}, where $\tilde{H}$ has the form \eqref{ech} or \eqref{nph}, satisfy \eqref{hs0}. Then, for $h_0<h<h_0+\beta$ we have
\begin{equation}
M(h,\delta)=(h-h_0)^{\frac{p+1}{2p}} \sum_{j\geq0}\tilde{b}_j(\delta) (h-h_0)^{\frac{j}p},
\end{equation}
where
$$\tilde{b}_j(\delta)=sgn(k)\Big(\dfrac{k}{D}\Big)^{\frac{1+2j-p}{2p}}b_j(\delta).$$

Especially, for p=1, namely system \eqref{e12} has an elementary center at the origin, then, for $h_0<h<h_0+\beta$ we have
\begin{equation*}
M(h,\delta)=(h-h_0) \sum_{j\geq0}\tilde{b}_j(\delta) (h-h_0)^j,
\end{equation*}
where
$$\tilde{b}_j(\delta)=sgn(k)\Big(\dfrac{k}{D}\Big)^j b_j(\delta).$$
}

As we known, $D\neq0$ and $k\neq0$. Then, by Theorem 3 and Theorem 4, to determine the number of limit cycles bifurcating from the center $(x_0,y_0)$ of system \eqref{e11}, we need only to consider the Melnikov function $\tilde{M}$ of near-Hamiltonian system \eqref{ths}.

Theorem 9 told us a method to study limit cycles bifurcating from a center $(x_0,y_0)$ of the near-Hamiltonian system \eqref{e11}. First, by introducing suitable linear change of variables or rescaling the time, we can transform system \eqref{e11} into a normal form \eqref{ths} near the center $(x_0,y_0)$. Then, using Theorem 2 after executing the Mathematica program, we can obtain the first several coefficients of the Melnikov function of system \eqref{ths} at the origin. Finally, by Theorems 3 and 4 we can find the number of limit cycles of near-Hamiltonian system \eqref{ths} at the origin, which is the same to that of the near-Hamiltonian system \eqref{e11} at the center $(x_0,y_0)$.

Suppose that for systems \eqref{e11} and \eqref{ths},
$$p_x+q_y=\sum_{i+j=0}^2c_{ij}x^iy^j,\;\tilde{p}_u+\tilde{q}_v=\sum_{i+j=0}^2\tilde{c}_{ij}u^iv^j.$$
Then, using the Mathematica program we can compute
\begin{equation}
\begin{array}{ll}
\tilde{c}_{00}&=\dfrac1k[c_{00}+x_0\ c_{10}+y_0\ c_{01}+x^2_0\ c_{20}+x_0y_0\ c_{11}+y^2_0\ c_{02}],\vspace{0.2cm}\\
\tilde{c}_{10}&=\dfrac1{kD}[-c\ (c_{01}+x_0\ c_{11}+2 y_0\ c_{02})+d\ (c_{10}+y_0\ c_{11}+2 x_0\ c_{20})],\vspace{0.2cm}\\
\tilde{c}_{01}&=\dfrac1{kD}[a\ (c_{01}+x_0\ c_{11}+2 y_0\ c_{02})-b\ (c_{10}+y_0\ c_{11}+2 x_0\ c_{20})],\vspace{0.2cm}\\
\tilde{c}_{20}&=\dfrac1{kD^2}[c^2\ c_{02}-cd\ c_{11}+d^2\ c_{20}],\vspace{00.2cm}\\
\tilde{c}_{11}&=\dfrac1{kD^2}[-2ac\ c_{02}+(bc+ad)\ c_{11}-2bd\ c_{20}],\vspace{00.2cm}\\
\tilde{c}_{02}&=\dfrac1{kD^2}[a^2\ c_{02}-ab\ c_{11}+b^2\ c_{20}].\vspace{00.2cm}\\
\end{array}
\end{equation}

Now, we turn to considering the near-Hamiltonian system obtained by perturbing system $(S)$ with cubic polynomials as follows:
\begin{equation}\label{ens}
\dot{x}=H_y +\varepsilon p(x,y),\quad \dot{y}=-H_x+\varepsilon q(x,y),
\end{equation}
where
\begin{equation*}
H(x,y)=\sum_{i+j=2,4} h_{ij} x^i y^j,
\end{equation*}
\begin{equation*}
p(x,y)=\sum_{i+j=1}^3 a_{ij}x^iy^j,\;q(x,y)=\sum_{i+j=1}^3 b_{ij}x^iy^j
\end{equation*}
and
$$p_x+q_y=c_{00}+c_{10}x+c_{01}y+c_{20}x^2+c_{11}xy+c_{20}y^2.$$

From section 3, if $\Delta=-4h^2_{11}-8h_{02}(h_{20}+h_{22})>0$, which means that $h_{02}\neq0$ and $h_{20}+h_{22}\neq0$, then there is an elementary center of Hamiltonian system $(S)$ at the point $(0,1)$. For convenience, we let $h_{20}+h_{22}=1/2$ and $\Delta=1$. Then, introducing a change of variables $(x,y)$ as
$$
u=-(y-1),\;v=x-2h_{11}(y-1),
$$
and then taking $(u,v)$ as $(x,y)$, one can obtain the normal form of system \eqref{ens} near the center $(0,1)$ as
\begin{equation}\label{hs1}
\dfrac{du}{dt}=\dfrac{\partial{H}_1}{\partial{v}}(u,v)+\varepsilon\tilde{p}_1(u,v),\;
\dfrac{dv}{dt}=-\dfrac{\partial{H}_1}{\partial{u}}(u,v)+\varepsilon\tilde{q}_1(u,v),
\end{equation}
where $H_1$ is the same as that in section 3 and
$$\dfrac{\partial\tilde{p}_1}{\partial{u}}+\dfrac{\partial\tilde{q}_1}{\partial{v}}=\sum_{i+j=0}^2\tilde{c}1_{ij}u^iv^j,$$
where
$$\begin{array}{ll}
\tilde{c}1_{00}&=c_{00}+c_{01}+c_{02},\;\tilde{c}1_{10}=-2 h_{11} c_{10}-c_{01}-2c_{02}-2h_{11}c_{11},\;
\tilde{c}1_{01}=c_{10}+c_{11},\vspace{0.1cm}\\
\tilde{c}1_{20}&=c_{02}+2h_{11}c_{11}+4h^2_{11}c_{20},\;\tilde{c}1_{11}=-c_{11}-4h_{11}c_{20},\;
\tilde{c}1_{02}=c_{20}.\vspace{0.01cm}
\end{array}$$
Obviously, $\tilde{c}1_{00},\tilde{c}1_{10},\tilde{c}1_{01},\tilde{c}1_{20},\tilde{c}1_{11},\tilde{c}1_{02}$ are all linear functions of $c_{00},c_{10},c_{01},c_{20},c_{11},c_{02}$ and are independent. Therefore, we can denote $\tilde{c}1_{ij}$ as $c_{ij}$ for $0\leq{i+j}\leq2,\,i\geq0,\,j\geq0$.\vspace{0.1cm}

For system \eqref{hs1}, there are too many parameters in the corresponding Hamiltonian function. Hence, the coefficients of the Melnikov function will be very complicated and it is very hard to give out the accurate number of limit cycles of system \eqref{hs1} near the origin. However, suppose that there are altogether $m$ parameters in $(\delta,a)$, by Theorem 3 or Theorem 4 we know that one can find at most $m-1$ limit cycles near the origin. As the result, one can find at most 9 limit cycles of system \eqref{hs1} near the origin. Here, we consider an simple case of system \eqref{hs1} to illustrate that system \eqref{hs1} can have at least 6 limit cycles near the origin. As an especial case, we assume that $h_{11}=h_{22}=0,\,h_{31}=1$ in system \eqref{hs1}. Namely, we consider the following system,
\begin{equation}\label{esc}
\dot{u}=H_v+\varepsilon{p}(u,v),\;\dot{v}=-H_u+\varepsilon{q}(u,v),
\end{equation}
where $p,q$ are the same as before and
$$
H(u,v)=\dfrac12 (u^2+v^2)-\dfrac12 u^3+v^3-\dfrac18 u^4- u v^3+ r v^4
$$
with $r=h_{40}$.

\vspace{0.10in} \noindent {\bf Theorem 7.} \ {\it For system
\eqref{esc}, introduce
$$
\delta=(c_{00},c_{10},c_{01},c_{02},c_{11},c_{20}), \ \sigma =(r),
$$
and
$$
\delta_0=(c_{00}^*,c_{10}^*,c_{01}^*,c_{02}^*,c_{11}^*,c_{20}^*), \
\sigma_0 =(r^*),
$$
where $\delta_0$ and $\sigma_0$ satisfy
$$
\begin{array}{c}
f(r^*)=0,\;c^*_{00}=0,\;c^*_{02}=-\dfrac32 c^*_{10}+3 c^*_{01}-c^*_{20},\;c^*_{01}\neq0,\vspace{0.2cm}\\
c^*_{11}=(12 r^*-\dfrac{191}4) c^*_{10}+(40 r^*+3) c^*_{01}+(8 r^*-\dfrac{47}2) c^*_{20},\;\vspace{0.2cm}\\
c^*_{10}=\dfrac{2}{ 17680 r^*-43347}\big[(1792 (r^*)^2- 17328 r^*-696 )c^*_{01}+(+11107  - 5680 r^*) c^*_{20}\big],\vspace{0.2cm}\\
c^*_{20}=\dfrac{6544683 +11682846 r^* +30378720 (r^*)^2 - 40278528 (r^*)^3 - 21471232 (r^*)^4}{2 (2770317 - 11417834 r^* + 11442688 (r^*)^2 + 442880 (r^*)^3)}c^*_{01}, \end{array}
$$
and
$$
f(r)=780900831 - 459924741 r - 2093583104 r^2 - 1597844992 r^3 - 108363776 r^4 + 1381957632 r^5.
$$
Then, for some $(\varepsilon,\delta,\sigma)$ near
$(0,\delta_0,\sigma_0)$ system \eqref{esc} has $6$ limit cycles near the origin.}

\vspace{0.10in} \noindent \emph{Proof.} \  Let $M(h,\delta)$ be the Melnikov function of system $(nH)_1$ near the origin. Then, by Theorem 2, we have
$$
M(h,\delta)=h\sum_{j\geq0}b_j(\delta)h^j.
$$
Executing the Mathematica
program yields
$$
b_0 = 2\, \pi \, c_{00}.
$$
Letting $\, c_{00}=0$, under which $\, b_0 = 0$, we
obtain the following expressions for $\, b_i$'s.
\begin{equation}\label{t61}
\begin{array}{rl}
b_{1}|_{b_0=0}=&\!\!\!\pi[c_{20}+c_{02}-3c_{01}+\dfrac32 c_{10}],\vspace{0.2cm}\\
b_{2}|_{b_0=0} =&\!\!\!\dfrac{\pi}{8}\big[(25 - 12 r)c_{10}+140 (-3 + 2 r)c_{01}+(25 - 4 r)c_{20}-4c_{11}+8 (18 - 5 r)c_{02}\big], \vspace{0.2cm}\\
b_i|_{b_0=0}=&\!\!\!\l_i(\delta_{i1}c_{10}+\delta_{i2}c_{01}+\delta_{i3}c_{20}+\delta_{i4}c_{11}+\delta_{i5}c_{02}),\ \ i=3,4,5,6
\end{array}
\end{equation}
where
$$
l_3=-\dfrac{5 \pi}{64},\;l_4=\dfrac{21 \pi}{512},\;l_5=-\dfrac{21 \pi}{2048},\;l_6=\dfrac{33 \pi}{8193},
$$
and
$$
\begin{array}{rl}
\delta_{31}=&\!\!\!939 - 210 r - 84 r^2,\;\delta_{32}=462(39 - 52 r + 12 r^2),\;\delta_{33}=-2(333 - 270 r + 28 r^2),\vspace{0.2cm}\\
\delta_{34}=&\!\!\!-2 (123 + 28 r),\;\delta_{35}=-8 (749 - 696 r + 63 r^2),\vspace{0.2cm}\\
\delta_{41}=&\!\!\!85233 - 85896 r + 6584 r^2 + 1056 r^3,\;\delta_{42}=3432 (-331 + 646 r - 340 r^2 + 40 r^3),\;\vspace{0.2cm}\\
\delta_{43}=&\!\!\!21278 - 34416 r + 13968 r^2 - 704 r^3,\;\delta_{44}=-8 (-3009 + 1706 r + 132 r^2),\;\vspace{0.2cm}\\
\delta_{45}=&\!\!\!64 (-5828 + 9101 r - 3222 r^2 + 143 r^3),\vspace{0.2cm}\\
\delta_{51}=&\!\!\!13788577 - 23251746 r + 9251880 r^2 - 382800 r^3 - 34320 r^4,\vspace{0.2cm}\\
\delta_{52}=&\!\!\!92378 (1821 - 4664 r + 3864 r^2 - 1120 r^3 + 80 r^4),\vspace{0.2cm}\\
\delta_{53}=&\!\!\!-2 (1156991 - 2579166 r + 1845528 r^2 - 391600 r^3 + 11440 r^4),\;\vspace{0.2cm}\\
\delta_{54}=&\!\!\!-2 (2028023 - 2589736 r + 582120 r^2 + 22880 r^3),\;\vspace{0.2cm}\\
\delta_{55}=&\!\!\!-32 (1717743 - 3736392 r + 2422830 r^2 - 462176 r^3 + 12155 r^4),\vspace{0.2cm}\\
\delta_{61}=&\!\!\!-1454668039 + 3393291252 r - 2432266980 r^2 + 541324640 r^3 -
 13999440 r^4 - 806208 r^5,\vspace{0.2cm}\\
\delta_{62}=&\!\!\!386308 (-44055 + 139650 r - 157944 r^2 + 75600 r^3 - 14000 r^4 +
   672 r^5),\;\vspace{0.2cm}\\
\delta_{63}=&\!\!\!2 (97965305 - 271364748 r + 269902620 r^2 - 108076960 r^3 +
   14155440 r^4 - 268736 r^5),\vspace{0.2cm}\\
\delta_{64}=&\!\!\!4 (109528617 - 212220190 r + 113651592 r^2 - 14093040 r^3 -
   335920 r^4),\;\vspace{0.2cm}\\
\delta_{65}=&\!\!\!128 (43337762 - 120692295 r + 114411600 r^2 - 42231670 r^3 +
   5071950 r^4 - 88179 r^5).\vspace{0.2cm}
\end{array}
$$

Let
\begin{equation}\label{t62}
\widetilde{b}_j=b_j|_{b_0=0},\ \ j=1,2,3,4,5,6.
\end{equation}
Noticing \eqref{t61}, we can solve $c_{02}, c_{11}$ from
$\widetilde{b}_1=\widetilde{b}_2=0$ as follows:
\begin{equation}\label{t63}
c_{02}=-\dfrac32 c_{10}+3 c_{01}-c_{20},\;
c_{11}=(12 r-\dfrac{191}4) c_{10}+(40 r+3) c_{01}+(8 r-\dfrac{47}2) c_{20}.
\end{equation}
Substituting \eqref{t63} into $\widetilde b_3,\,\widetilde b_4\, \widetilde b_5\,$
and $\, \widetilde b_6\,$ results in
\begin{equation}\label{t64}
\begin{array}{rl}
\widetilde b_{3} =&\!\!\!\dfrac{5 \pi}{128} \big[(17680 r-43347 )c_{10}+16(87 + 2166 r - 224 r^2)c_{01}+2(-11107 + 5680r)c_{20}\big],\vspace{0.2cm}\\
\widetilde b_{i} =&\!\!\! \dfrac{l_i}{2}(\delta'_{i1}c_{10}+\delta'_{i2}c_{01}+\delta'_{i3}c_{20}),\;i=4,5,6
\end{array}\end{equation}
where $l_i$ same as before and
$$
\begin{array}{rl}
\delta'_{41}=&\!\!\!-2 (1794159 - 1900148 r + 429248 r^2),\;\delta'_{42}=64 (1725 + 43488 r - 34292 r^2 + 2112 r^3),\vspace{0.2cm}\\
\delta'_{43}=&\!\!\!-4 (458703 - 530676 r + 138304 r^2),\;\vspace{0.2cm}\\
\delta'_{51}=&\!\!\!579832875 - 997181804 r + 486587688 r^2 - 68706176 r^3,\vspace{0.2cm}\\
\delta'_{52}=&\!\!\!16 (-1106391 - 27357598 r + 41005384 r^2 - 13225168 r^3 + 549120 r^4),\vspace{0.2cm}\\
\delta'_{53}=&\!\!\!-2 (-147970875 + 268572172 r - 142634920 r^2 + 22244992 r^3),\;\vspace{0.2cm}\\
\delta'_{61}=&\!\!\!-4 (15347742095 - 36178820899 r + 28146658686 r^2 - 8398426168 r^3 +
   800059520 r^4),\vspace{0.2cm}\\
\delta'_{62}=&\!\!\!32 (58577817 + 1411250460 r - 3104524228 r^2 + 1937691360 r^3 -
   357475040 r^4 + 10749440 r^5),\;\vspace{0.2cm}\\
\delta'_{63}=&\!\!\!-8 (3911748231 - 9589874467 r + 7894793822 r^2 - 2537774136 r^3 +
   260074880 r^4).
\end{array}
$$

For $r$ near $r^*$ we have $17680 r-43347\neq 0$. Thus, $\widetilde{b}_3=0$ if
and only if
$$
c_{10}=\dfrac{2}{ 17680 r-43347}\big[(1792 r^2- 17328 r-696 )c_{01}+(+11107  - 5680 r) c_{20}\big],
$$
which results in
$$
\begin{array}{ll}
\widetilde{b}_4=&\dfrac{21 \pi}{32 (-43347 + 17680 r)}\big[(6544683 + 11682846 r+ 30378720 r^2 - 40278528 r^3\vspace{0.2cm}\\
&- 21471232 r^4)c_{01}-2(2770317 - 11417834 r + 11442688 r^2 + 442880 r^3)c_{20}\big].
\end{array}
$$

For $r$ near $r^*$ we have $2770317 - 11417834 r + 11442688 r^2 + 442880 r^3\neq 0$. Thus, $\widetilde{b}_4=0$ if
and only if
$$
c_{20}=\dfrac{6544683 +11682846 r +30378720 r^2 - 40278528 r^3 - 21471232 r^4}{2 (2770317 - 11417834 r + 11442688 r^2 + 442880 r^3)}c_{01},
$$
which results in
$$
\widetilde{b}_{5+i}=\Delta_i L,\;i=0,1,
$$
where
$$
\begin{array}{rl}
L=&\!\!\!\dfrac{3 \pi (4 r^2+ 4 r-7) c_{01}}{32 (2770317 - 11417834 r + 11442688 r^2 + 442880 r^3)},\vspace{0.2cm}\\
\Delta_0=&\!\!\!-14 f(r),\vspace{0.2cm}\\
\Delta_1=&\!\!\!11 (-77526899979 + 74770117515 r + 170378796326 r^2 +   50680931328 r^3\vspace{0.2cm}\\
&\!\!\! - 71020059648 r^4- 138129047552 r^5 +   61673963520 r^6).
\end{array}
$$

Obviously, there is at least one real solution to equation $f(r)=0$ and denote it by $r^*$. Then, for $r$ near $r^*$ and $c_{01}$ near $c^*_{01}$ we have $L\neq0$. Thus, $\widetilde{b}_5=0$ if and only if $\Delta_0=0$. When $r=r^*$, we have $\widetilde{b}_5=0$ but $\Delta_1\neq0$, which results in $\widetilde{b}_6\neq0$.
Hence, it is clear that the conclusion follows from Theorem 4.

The proof is complete.
\vspace{0.5cm}

Therefore, if system $(S)$ has an elementary center at $(0,1)$ (or $(0,-1)$), then one can find 6 to 9 limit cycles of system \eqref{ens} near the center $(0,1)$ (or $(0,-1)$).\vspace{0.3cm}

If $\triangle=-4h^2_{11}-8h_{02}(h_{20}+h_{22})=0$, $h_{31}=-2h_{11}h_{20}$ and $h^2_{20}+h_{40}<0$, then there is a nilpotent center of order 1 of Hamiltonian system $(S)$ at the point $(0,1)$. For convenience, let $h_{02}=1/2$. Then, we can take
$$h_{22}=-h^2_{11}-h_{20}.$$
Introducing a change of variables $(x,y)$ as
$$
u=-\frac12\ (y-1),\;v=h_{11}x+y-1,
$$
one can obtain the normal form of system \eqref{ens} near the center $(0,1)$ as
\begin{equation}\label{ns2}
\dfrac{du}{dt}=\dfrac{\partial{H}_2}{\partial{v}}(u,v)+\varepsilon\tilde{p}_2(u,v),\;
\dfrac{dv}{dt}=-\dfrac{\partial{H}_2}{\partial{u}}(u,v)+\varepsilon\tilde{q}_2(u,v),
\end{equation}
where $H_2$ is the same as that in section 3, namely,
\begin{equation*}
H_2(u,v)=\dfrac12v^2+2(2h_{20}-h^2_{11})u^2v+\dfrac12v^3+(2h^4_{11}-8h^2_{11}h_{20}-8h_{40})u^4+(2h_{20}-h^2_{11})u^2v^2+\dfrac18v^4,
\end{equation*}
and
$$\dfrac{\partial\tilde{p}_2}{\partial{u}}+\dfrac{\partial\tilde{q}_2}{\partial{v}}=\sum_{i+j=0}^2\tilde{c}2_{ij}u^iv^j,$$
where
$$
\begin{array}{ll}
\tilde{c}2_{00}&=c_{00}+c_{01}+c_{02},\;\tilde{c}2_{10}=-2 c_{10}+2 h_{11}c_{01}-2c_{11}+4h_{11}c_{02},\;
\tilde{c}2_{01}=c_{01}+2c_{02},\vspace{0.2cm}\\
\tilde{c}2_{20}&=4c_{20}-4h_{11}c_{11}+4h^2_{11}c_{02},\;\tilde{c}2_{11}=-2c_{11}+4h_{11}c_{02},\;
\tilde{c}2_{02}=c_{02}.
\end{array}
$$
Obviously, $\tilde{c}2_{00},\tilde{c}2_{10},\tilde{c}2_{01},\tilde{c}2_{20},\tilde{c}2_{11},\tilde{c}2_{02}$ are all linear functions of $(c_{00},c_{10},c_{01},c_{20},c_{11},c_{02})$ and are independent. Therefore, we can denote $\tilde{c}2_{ij}$ as $c_{ij}$ for $0\leq{i+j}\leq2,\,i\geq0,\,j\geq0$.

\vspace{0.10in} \noindent {\bf Theorem 8.} \ {\it For system
\eqref{ns2}, we introduce
$$
\delta=(c_{00},c_{10},c_{01},c_{02},c_{11},c_{20}), \ \sigma =(A,B),
$$
and
$$
\delta_0=(c_{00}^*,c_{10}^*,c_{01}^*,c_{02}^*,c_{11}^*,c_{20}^*), \
\sigma_0 =(A^*,B^*),
$$
where $\delta_0$ and $\sigma_0$ satisfy
$$
 A^*=0, \
B^*>0,  \ c_{00}^*=0 ,\,c^*_{20}=0, \
c_{02}^*=\dfrac{3}{2}c_{01}^*,\;c_{01}^*\neq0.
$$
Then, for some $(\varepsilon,\delta,\sigma)$ near
$(0,\delta_0,\sigma_0)$ the system \eqref{ns2} has $4$ limit cycles near the origin.}

\vspace{0.10in} \noindent \emph{Proof.} \  Let $M(h,\delta)$ be the Melnikov function of system \eqref{ns2} near the origin. Then, by Theorem 2, we have
$$
M(h,\delta)=h^{3/4}\sum_{j\geq0}b_j(\delta)h^{j/2}.
$$
Executing the Mathematica
program yields
$$
b_0 =\dfrac{2\Gamma(5/4)\sqrt{2 \pi}}{\Gamma(7/4)(B-2A^2)^{1/4}}c_{00}.
$$
Letting $c_{00}=0$, under which $\, b_0 = 0$, we
obtain the following expressions for $\, b_i$'s.
\begin{equation}\label{b123}
\begin{array}{rl}
b_{1}|_{b_0=0}=&\!\!\!\dfrac{\Gamma(3/4)\sqrt{\pi/2}}{\Gamma(9/4)(B-2A^2)^{3/4}}(c_{20}-2Ac_{01}),\vspace{0.2cm}\\
b_{2}|_{b_0=0} =&\!\!\!\dfrac{\Gamma(1/4)\sqrt{\pi/2}}{2\Gamma(11/4)(B-2A^2)^{5/4}}[(A^2-B)(3c_{01}-2c_{02})+Ac_{20}],\vspace{0.2cm}\\
b_3|_{b_0=0}=&\!\!\!\dfrac{\Gamma(7/4)\sqrt{2\pi}}{\Gamma(13/4)(B-2A^2)^{7/4}}[2A(A^2-B)(5c_{01}-4c_{02})+(A^2+B)c_{20}],\vspace{0.2cm}\\
b_4|_{b_0=0}=&\!\!\!-\dfrac{5\Gamma(1/4)\sqrt{2\pi}}{16\Gamma(15/4)(B-2A^2)^{9/4}}[(5A^4-6A^2B+3B^2)(7c_{01}-6c_{02})+2A(A^2-3B)c_{20}].
\end{array}
\end{equation}

Let
\begin{equation}\label{e36}
\widetilde{b}_j=b_j|_{b_0=0},\ \ j=1,2,3.
\end{equation}
Noticing \eqref{b123}, we can obtain $c_{20}=2Ac_{01}$ by solving $\widetilde{b}_1=0$.
Substituting $c_{20}=2Ac_{01}$ into $\widetilde b_2,\,\widetilde b_3\, $ and $\widetilde b_4$
results in
$$\begin{array}{rl}
\widetilde{b}_2=&\!\!\! \dfrac{\Gamma(5/4)\sqrt{2\pi}}{\Gamma(11/4)(B-2A^2)^{5/4}}[(5A^2-3B)c_{01}-2(A^2-B)c_{02}],\vspace{0.2cm}\\
\widetilde{b}_3
=&\!\!\!\dfrac{3A\Gamma(3/4)\sqrt{2\pi}}{\Gamma(13/4)(B-2A^2)^{7/4}}[(3A^2-2B)c_{01}-2(A^2-B)c_{02}],\vspace{0.2cm}\\
\widetilde{b}_4=&\!\!\!-\dfrac{15\Gamma(1/4)\sqrt{2\pi}}{16\Gamma(15/4)(B-2A^2)^{9/4}}[(13A^4-18A^2B+7B^2)c_{01}-2(5A^4-6A^2B+3B^2)c_{02})].\vspace{0.2cm}
\end{array}
$$

For $(A,B)$ near $(A^*,B^*)$, we have $A^2-B<0$. Then $\widetilde{b}_2=0$ if and only if
$$c_{02}=\dfrac{5A^2-3B}{2(A^2-B)}c_{01},$$
which yields
$$
\begin{array}{rl}
\widetilde{b}_3
=&\!\!\!\dfrac{3A\Gamma(3/4)\sqrt{2\pi}}{\Gamma(13/4)(B-2A^2)^{3/4}}c_{01},\vspace{0.2cm}\\
\widetilde{b}_4=&\!\!\!-\dfrac{6\Gamma(9/4)\sqrt{2\pi}}{\Gamma(15/4)(A^2-B)(B-2A^2)^{5/4}}(3A^4-2A^2B+B^2)c_{01}.\vspace{0.2cm}
\end{array}
$$
Thus, it follows that
when $\widetilde{b}_3=0$ we have $A=0$ or $c_{01}=0$. However, $c_{01}$ is a common factor of $\widetilde{b}_3$ and $\widetilde{b}_4$. Taking $A=A^*=0$, we have
$$\widetilde{b}_4=\,\,\dfrac{6\Gamma(9/4)\sqrt{2\pi}}{\Gamma(15/4)B^{1/4}}c_{01}.$$
Hence, it is clear that the conclusion follows from Theorem 4.

The proof is complete.
\vspace{0.3cm}

Therefore, if system $(S)$ has a nilpotent center of order 1 at $(0,1)$ (or $(0,-1)$), then one can find 4 limit cycles of system \eqref{ens} near the center $(0,1)$ (or $(0,-1)$).

\vspace{0.10in}
\section{Bifurcation in the symmetric near-Hamiltonian systems}

In this section, we will consider limit cycles bifurcating from the symmetric centers of some symmetric near-Hamiltonian systems. We said a near-hamiltonian system \eqref{e11} is symmetric, which means that not only the corresponding Hamiltonian function is symmetric, but also the perturbing terms $p(x,y)$ and $q(x,y)$ are also symmetric with respect to the origin.

Now, we consider the following symmetric near-Hamiltonian system:
\begin{equation}\label{sns}
\dot{x}=H_y +\varepsilon p(x,y),\quad \dot{y}=-H_x+\varepsilon q(x,y),
\end{equation}
where
\begin{equation*}
H(x,y)=\sum_{i+j=2,4} h_{ij} x^i y^j,
\end{equation*}

\begin{equation*}
p(x,y)=\sum_{i+j=1,3}a_{ij}x^iy^j,\;q(x,y)=\sum_{i+j=1,3}b_{ij}x^iy^j
\end{equation*}
and
$$p_x+q_y=c_{00}+c_{20}x^2+c_{11}xy+c_{20}y^2.$$

Obviously, system \eqref{sns} is obtained by perturbing symmetric Hamiltonian system $(S)$ with symmetric cubic polynomials and is symmetric with respect to the origin. If there are $m$ limit cycles bifurcating from the center $(x_0,y_0)$ (not the origin) of the corresponding Hamiltonian system, then it is same to the symmetric center $(-x_0,-y_0)$. Namely, there are $2m$ limit cycles of the symmetric near-Hamiltonian system.\vspace{0.1cm}

From section 3, if $\Delta=-4h^2_{11}-8h_{02}(h_{20}+h_{22})>0$, which means that $h_{02}\neq0$ and $h_{20}+h_{22}\neq0$, then there is an elementary center of Hamiltonian system $(S)$ at the point $(0,1)$. For convenience, we let $h_{20}+h_{22}=1/2$ and $\Delta=1$. Then, introducing a change of variables $(x,y)$ as
$$
u=-(y-1),\;v=x-2h_{11}(y-1),
$$
one can obtain the normal form of system \eqref{sns} near the center $(0,1)$ as
\begin{equation}\label{sns1}
\dfrac{du}{dt}=\dfrac{\partial{H}_1}{\partial{v}}(u,v)+\varepsilon\hat{p}_1(u,v),\;
\dfrac{dv}{dt}=-\dfrac{\partial{H}_1}{\partial{u}}(u,v)+\varepsilon\hat{q}_1(u,v),
\end{equation}
where $H_1$ is the same as that in section 3 and
$$\dfrac{\partial\hat{p}_1}{\partial{u}}+\dfrac{\partial\hat{q}_1}{\partial{v}}=\sum_{i+j=0}^2\hat{c}1_{ij}u^iv^j,$$
where
$$
\begin{array}{ll}
\hat{c}1_{00}&=c_{00}+c_{02},\;\hat{c}1_{10}=-2c_{02}-2h_{11}c_{11},\;
\hat{c}1_{01}=c_{11},\vspace{0.2cm}\\
\hat{c}1_{20}&=c_{02}+2h_{11}c_{11}+4h^2_{11}c_{20},\;\hat{c}1_{11}=-c_{11}-4h_{11}c_{20},\;
\hat{c}1_{02}=c_{20}\vspace{0.2cm}
\end{array}
$$

For simplification, we study system \eqref{sns1} in two cases $h_{11}=0$ and $h_{11}=1$ separately. First, let $h_{11}=0$. Then for system \eqref{sns1} we have
$$H_1(u,v)=\dfrac12(u^2+v^2)-\dfrac12 u^3-2 h_{22} u v^2+h_{31} v^3+\dfrac18 u^4+h_{22}u^2 v^2-h_{31} u v^3+h_{40} v^4$$
and
$$
\begin{array}{ll}
\hat{c}1_{00}&=c_{00}+c_{02},\;\hat{c}1_{10}=-2c_{02},\;
\hat{c}1_{01}=c_{11},\vspace{0.2cm}\\
\hat{c}1_{20}&=c_{02},\;\hat{c}1_{11}=-c_{11},\;
\hat{c}1_{02}=c_{20}.\vspace{0.2cm}
\end{array}
$$

Therefore, we have the following theorem

\vspace{0.10in} \noindent {\bf Theorem 9.} \ {\it Suppose system
\eqref{sns1} satisfies $h_{11}=0$ and
$$h_{02}(h_{20}+h_{22})<0.$$
Introduce
$$
\delta=(c_{00},c_{02},c_{11},c_{20}), \ \sigma =(h_{31},h_{22},h_{40}),
$$
and
$$
\delta_0=(c_{00}^*,c_{02}^*,c_{11}^*,c_{20}^*), \
\sigma_0 =(h_{31}^*,h_{22}^*,h_{40}^*),
$$
where $\delta_0$ and $\sigma_0$ satisfy
$$
\begin{array}{c}c_{00}^*=-c_{02}^*,\ c_{20}^*=3 h_{31}^* c_{11}^* +(2+ 4 h_{22}^*) c_{02}^*,\ c_{11}^*\neq0,\ \vspace{0.25cm}\\
c_{02}^*=\dfrac{h_{31}^*[20(h_{22}^*)^2-8 h_{22}^*-10 h_{40}^*-1]}{-2+20 (h_{31}^*)^2+8h_{22}^*-4h_{40}^*}c_{11}^*,\vspace{0.15cm}\\ [1.5ex]
h_{31}^*\neq0, \ h_{22}^*=\dfrac12-\dfrac{105\pm56\sqrt3}{4}(h_{31}^*)^2,  \ h_{40}^*=\dfrac1{14}[28(h_{22}^*)^2-21(h_{31}^*)^2-12 h_{22}^*-1].
 \end{array}
$$
Then, for some $(\varepsilon,\delta,\sigma)$ near
$(0,\delta_0,\sigma_0)$ the system \eqref{sns1} has $5$ limit cycles near the origin.}

\vspace{0.10in} \noindent \emph{Proof.} \  Let $M(h,\delta)$ be the Melnikov function of system \eqref{sns1} near the origin. Then, by Theorem 2, we have
$$
M(h,\delta)=h\sum_{j\geq0}b_j(\delta)h^j.
$$
Executing the Mathematica program yields
$$
b_0 = 2\, \pi \, (c_{00}+c_{02}).
$$
Letting $c_{00}=-c_{02}$, which implies $\, b_0 = 0$, we
obtain the following expressions for $\, b_i$'s:
\begin{equation}\label{e135}
\begin{array}{rl}
b_{1}|_{b_0=0}&\!\!\!=\pi[c_{20}-3 h_{31} c_{11}-(2+4h_{22})c_{02}],\vspace{0.2cm}\\
b_{2}|_{b_0=0}&\!\!\! =\dfrac{\pi}{2}\big[(1 + 4 h_{22} + 20 h_{22}^2 + 35 h_{31}^2 - 10 h_{40}) c_{20}\vspace{0.2cm}\\
&\!\!\!
+(h_{31} + 20 h_{22} h_{31} - 140 h_{22}^2h_{31} - 105 h_{31}^3+ 70h_{31}h_{40})c_{11} \vspace{0.2cm}\\
&\!\!\!+(-10 + 10 h_{31}^2 - 4 h_{22} (3 + 6 h_{22} + 20 h_{22}^2 + 35 h_{31}^2 - 10 h_{40}) + 4 h_{40})c_{02}\big],\vspace{0.2cm}\\
b_i|_{b_0=0}&\!\!\!=L_i(\delta_{i1}c_{20}+\delta_{i2}c_{11}+\delta_{i3}c_{02}),\ \ i=3,4,5,
\end{array}
\end{equation}
where
$$
L_3=-\dfrac{5}{32} \pi,\;L_4=-\dfrac{7}{64} \pi,\;L_5=-\dfrac{21}{1024} \pi,
$$
and
$$
\begin{array}{rl}
\delta_{31}=&\!\!\!-7 - 8 h_{22} (3 + h_{22} (9 + 14 h_{22} (2 + 9 h_{22}))) - 3003 h_{31}^4 \vspace{0.2cm}\\
&\!\!\!+ 12 h_{40} + 112 h_{22} (1 + 9 h_{22}) h_{40} - 252 h_{40}^2  \vspace{0.2cm}\\
&\!\!\!+ 14 h_{31}^2 (1 + 36 (1 - 11 h_{22}) h_{22} + 198 h_{40})\vspace{0.2cm}\\
\delta_{32}=&\!\!\!h_{31} (-3 - 24 h_{22} + 28 h_{40} + 21 (-96 h_{22}^3 + 528 h_{22}^4 + 429 h_{31}^4,\vspace{0.2cm}\\
&\!\!\! + 2 h_{31}^2 (3 + 44 h_{22} (-3 + 13 h_{22}) - 286 h_{40}) +
    48 h_{22} h_{40} \vspace{0.2cm}\\
&\!\!\!+ 132h_{40}^2 - 8h_{22}^2(1 + 66 h_{40}))),\vspace{0.2cm}\\
\delta_{33}=&\!\!\!-2 (-63 - 2 h_{22} (35 + 4 h_{22}(15 + 2 h_{22} (15 + 7 h_{22} (5 + 18 h_{22})))) \vspace{0.2cm}\\
&\!\!\!+ 231 (3 - 26 h_{22}) h_{31}^4 + 12 h_{40} +
  24 h_{22} (3 + 14 h_{22} (1 + 6 h_{22})) h_{40}\vspace{0.2cm}\\
&\!\!\! - 28 (1 + 18 h_{22}) h_{40}^2 -
  6 h_{31}^2 (-1 + 14 h_{22} (-1 + 6 h_{22} (-3 + 22 h_{22}) \vspace{0.2cm}\\
&\!\!\!- 66 h_{40}) + 42 h_{40})),\vspace{0.2cm}\\
\delta_{41}=&\!\!\!\cdots.\vspace{0.2cm}
\end{array}
$$

Let
\begin{equation}\label{e136}
\widetilde{b}_j=b_j|_{b_0=0},\ \ j=1,2,3,4,5.
\end{equation}
Noticing \eqref{e135}, we can solve $c_{20}$ from
$\widetilde{b}_1=0$ as
\begin{equation}\label{e137}
c_{20}=3 h_{31} c_{11}+(2+4h_{22})c_{02}
\end{equation}
Substituting \eqref{e137} into $\widetilde b_2$ results in
$$
\widetilde{b}_2=2\pi\big[2(10 h_{31}^2+4 h_{22}^2-2 h_{40}-1)c_{02}+h_{31}(10 h_{40}+8h_{22}-20h_{22}^2+1)c_{11}\big].
$$
For $\sigma$ near $\sigma^*$ we have $10 h_{31}^2+4 h_{22}^2-2 h_{40}-1\neq 0$. Thus, we see that $\widetilde{b}_2=0$ if and only if
\begin{equation}\label{e138}
c_{02}=\dfrac{h_{31}(20h_{22}^2-8 h_{22}-10 h_{40}-1)}{-2+20 h_{31}^2+8h_{22}-4h_{40}}c_{11}.
\end{equation}

Substituting \eqref{e137} and \eqref{e138} into $\widetilde{b}_3,\,\widetilde{b}_4$ and $\widetilde{b}_5$ results in
\begin{equation}\label{e139}
\begin{array}{l}\widetilde b_{3+i} \equiv \Delta_{i} L, \ \
i=0,1,2,
\end{array}
\end{equation}
where
$$L=\dfrac{h_{31}(1-8h_{31}^2+16h_{31}^2h_{22}-8h_{22}^2+16h_{22}^4+4h_{40}-16h_{22}^2h_{40}+4h_{40}^2)}
{16(10 h_{31}^2+4h_{22}^2-2 h_{40}-1)}\pi c_{11},$$
and
$$\begin{array}{rl}
\Delta_0=&\!\!\! 80(1 + 12h_{22} - 28 h_{22}^2 + 21 h_{31}^2 + 14 h_{40}),\vspace{0.2cm}\\
\Delta_1=&\!\!\!168(3 + 32 h_{22} - 40 h_{22}^2 + 192 h_{22}^3 - 528 h_{22}^4 + 36 h_{31}^2\vspace{0.2cm}\\
&\!\!\! + 264 h_{22} h_{31}^2 + 429 h_{31}^4 + 28 h_{40} - 96 h_{22} h_{40} \vspace{0.2cm}\\
&\!\!\! + 528 h_{22}^2 h_{40} - 132 h_{40}^2),\vspace{0.2cm}\\
\Delta_2=&\!\!\!21 (143 + 1452 h_{22} - 1340 h_{22}^2 + 9632 h_{22}^3 - 16368 h_{22}^4 + 45760 h_{22}^5 \vspace{0.2cm}\\
&\!\!\! - 137280 h_{22}^6 + 1509 h_{31}^2 + 10648 h_{22} h_{31}^2 -  17160 h_{22}^2 h_{31}^2  \vspace{0.2cm}\\
&\!\!\!+ 205920 h_{22}^3 h_{31}^2 - 194480 h_{22}^4 h_{31}^2 + 11869 h_{31}^4 + 48620 h_{22} h_{31}^4  \vspace{0.2cm}\\
&\!\!\!+ 184756 h_{22}^2 h_{31}^4 +  138567 h_{31}^6 + 1158 h_{40} - 4528 h_{22} h_{40}  \vspace{0.2cm}\\
&\!\!\!+ 18832 h_{22}^2 h_{40} - 45760 h_{22}^3 h_{40} + 205920 h_{22}^4 h_{40} \vspace{0.2cm}\\
&\!\!\! + 2860 h_{31}^2 h_{40} -  102960 h_{22} h_{31}^2 h_{40} + 194480 h_{22}^2 h_{31}^2 h_{40} \vspace{0.2cm}\\
&\!\!\! - 92378 h_{31}^4 h_{40} - 5324 h_{40}^2 + 11440 h_{22} h_{40}^2 - \vspace{0.2cm}\\
&\!\!\!  102960 h_{22}^2 h_{40}^2 - 48620 h_{31}^2 h_{40}^2 + 17160 h_{40}^3). \vspace{0.15cm}
\end{array}
$$

For $c_{11}$ near $c_{11}^*$ and $\sigma$ near $\sigma^*$ we have $L\neq 0$. Thus, $\widetilde{b}_3=0$ if
and only if $\Delta_0=0$, which yields
$$h_{40}=-\dfrac1{14}(1 + 12h_{22} - 28 h_{22}^2 + 21 h_{31}^2).$$
In this case, we have
$$
\Delta_1=\dfrac{96}7\big[4 (1 - 2 h_{22})^2 + 420 (-1 + 2 h_{22}) h_{31}^2 + 1617 h_{31}^4\big].
$$
Then, $\Delta_1=0$ if and only if $h_{22}=f_1(h_{31})$ or $h_{22}=f_2(h_{31})$, where
$$f_i(h_{31})=\dfrac12-\dfrac{105}{4}h_{31}^2-(-1)^i14\sqrt3\, h_{31}^2,\;i=1,2.$$
Thus, it follows that when $h_{22}=f_i(h_{31})$ we have
$$
\Delta_2=688128 (54 + (-1)^i 31 \sqrt{3}) h_{31}^6,\;i=1,2.
$$
Hence, it is clear that the conclusion follows from Theorem 4.

The proof is complete.
\vspace{0.3cm}

Now we turn to the case $h_{11}=1$. Similarly, we can prove the following theorem

\vspace{0.10in} \noindent {\bf Theorem 9$'$.} \ {\it Suppose system
\eqref{sns1} satisfies $h_{11}=1$ and
$$-4-8h_{02}(h_{20}+h_{22})>0.$$
Introduce
$$
\delta=(c_{00},c_{02},c_{11},c_{20}), \ \sigma =(h_{31},h_{22},h_{40}),
$$
and
$$
\delta_0=(c_{00}^*,c_{02}^*,c_{11}^*,c_{20}^*), \
\sigma_0 =(h_{31}^*,h_{22}^*,h_{40}^*),
$$
where $\delta_0$ and $\sigma_0$ satisfy
$$
\begin{array}{c} h_{22}^*<\dfrac{7}{10}, \
h_{31}^*=\dfrac8{25}-\dfrac{4}{5}h_{22}^*-\dfrac1{25}\sqrt{\dfrac{1}{232}}\sqrt{(15-8\sqrt3)[7-10 (h_{22}^*)^2]},\vspace{0.15cm}\\
\ h_{40}^*=-\dfrac{139}{250}-\dfrac{22}{5}(h_{22}^*)^2+\dfrac{554}{175}h_{22}^*-12 h_{22}^* h_{31}^*+\dfrac{22}5h_{31}^*-\dfrac{17}2(h_{13}^*)^2,\vspace{0.15cm}\\ [1.5ex]
c_{00}^*=-c_{02}^*,\ c_{20}^*=\dfrac15\big[(52 h_{22}^*+60 h_{31}^*-22) c_{02}^* +(75 h_{31}^*+ 60 h_{22}^*-26) c_{11}^*\big],\
c_{02}^*\neq0,\vspace{0.15cm}\\ [1.5ex]
c_{11}^*=\dfrac{-2 L^* c_{02}^*}{(-8 + 20 h_{22}^* + 25 h_{31}^*) (-263 + 20 (42 - 25 h_{22}^*) h_{22}^* + 500 h_{31}^* + 1250h_{40}^*)} \end{array}
$$
with
$$
\begin{array}{ll}
L^*=&979 + 6250 (h_{31}^*)^2 - 50h_{31}^* (101 + 4h_{22}^*(-72 + 25h_{22}^*) - 250h_{40}^*) - 4050h_{40}^*\vspace{0.2cm} \\
&+ 4h_{22}^*(-1366 + 5*(457 - 200h_{22}^*)h_{22}^* + 2500h_{40}^*).
\end{array}
$$
Then, for some $(\varepsilon,\delta,\sigma)$ near
$(0,\delta_0,\sigma_0)$ the system \eqref{sns1} has $5$ limit cycles near the origin.}

\vspace{0.3cm}

Therefore, system \eqref{sns} has $5$ limit cycles near the elementary center $(0,1)$. For the symmetry, system \eqref{sns} also has 5 limit cycles near the center $(0,-1)$. Namely, there are altogether 10 limit cycles of system \eqref{sns} bifurcating from both symmetric elementary centers $(0,\pm1)$.

\vspace{0.3cm}

If $\triangle=-4h^2_{11}-8h_{02}(h_{20}+h_{22})=0$, $h_{31}=-2h_{11}h_{20}$ and $h^2_{20}+h_{40}<0$, then there is a nilpotent center of order 1 of Hamiltonian system $(S)$ at the point $(0,1)$. For convenience, let $h_{02}=1/2$. Then, we can take
$$h_{22}=-h^2_{11}-h_{20}.$$
Introducing a change of variables $(x,y)$ as
$$
u=-\frac12\ (y-1),\;v=h_{11}x+y-1,
$$
one can obtain the normal form of system \eqref{sns} near the center $(0,1)$ as
\begin{equation}\label{sns2}
\dfrac{du}{dt}=\dfrac{\partial{H}_2}{\partial{v}}(u,v)+\varepsilon\hat{p}_2(u,v),\;
\dfrac{dv}{dt}=-\dfrac{\partial{H}_2}{\partial{u}}(u,v)+\varepsilon\hat{q}_2(u,v),
\end{equation}
where $H_2$ is the same as that in section 3, namely,
\begin{equation*}
H_2(u,v)=\dfrac12 v^2+(4h_{20}-2h^2_{11}) u^2 v+\dfrac12 v^3+2(h^4_{11}-4h^2_{11}h_{20}-4h_{40}) u^4+(2h_{20}-h^2_{11}) u^2 v^2+\dfrac18 v^4,
\end{equation*}
and
$$\dfrac{\partial\hat{p}_2}{\partial{u}}+\dfrac{\partial\hat{q}_2}{\partial{v}}=\sum_{i+j=0}^2\hat{c}2_{ij}u^iv^j,$$
where
$$
\begin{array}{ll}
\hat{c}2_{00}&=c_{00}+c_{02},\;\hat{c}2_{10}=-2c_{11}+4h_{11}c_{02},\;
\hat{c}2_{01}=2c_{02},\vspace{0.2cm}\\
\hat{c}2_{20}&=4c_{20}-4h_{11}c_{11}+4h^2_{11}c_{02},\;\hat{c}2_{11}=-2c_{11}+4h_{11}c_{02},\;
\hat{c}2_{02}=c_{02}.
\end{array}
$$

\vspace{0.10in} \noindent {\bf Theorem 10.} \ {\it Suppose system
\eqref{sns2} satisfies
$$-4h^2_{11}-8h_{02}(h_{20}+h_{22})=0,\;h_{31}=-2h_{11}h_{20},\;h^2_{20}+h_{40}<0.$$
Introduce
$$
\delta=(c_{00},c_{02},c_{11},c_{20}), \ \sigma =(h_{11},h_{20},h_{40}),
$$
and
$$
\delta_0=(c_{00}^*,c_{02}^*,c_{11}^*,c_{20}^*), \
\sigma_0 =(h^*_{11},h^*_{20},h^*_{40}),
$$
where $\delta_0$ and $\sigma_0$ satisfy
$$
\begin{array}{c}
(h^*_{20})^2+h^*_{40}<0,  \ c_{00}^*=-c_{02}^*,\;c_{02}^*\neq0,\\ [1.5ex]
c^*_{20}=h^*_{11}c^*_{11}-2c^*_{02}[(h^*_{11})^2-h^*_{20}]. \end{array}
$$
Then, for some $(\varepsilon,\delta,\sigma)$ near
$(0,\delta_0,\sigma_0)$ the system \eqref{sns2} has $2$ limit cycles near the origin.}

\vspace{0.10in} \noindent \emph{Proof.} \ Let $M(h,\delta)$ be the Melnikov function of system \eqref{sns2} near the origin. Then, by Theorem 2, we have
$$
M(h,\delta)=h^{3/4}\sum_{j\geq0}b_j(\delta)h^{j/2}.
$$ Executing the Mathematica
program yields
$$
b_0 =\dfrac{2\Gamma(5/4)\sqrt{2 \pi}}{\Gamma(7/4)(B-2A^2)^{1/4}}(c_{00}+c_{02}).
$$
Letting $c_{00}=-c_{02}$ (implying $\, b_0 = 0$), we
obtain the following expressions for $\, b_i$'s:
\begin{equation}\label{b1111}
\begin{array}{rl}
b_{1}|_{b_0=0}=&\!\!\!\dfrac{\Gamma(3/4)\sqrt{\pi}}{\Gamma(9/4)(-2 h^2_{20}-2h_{40})^{3/4}}(c_{20}-h_{11}c_{11}+2c_{02}(h^2_{11}-h_{20})),\vspace{0.2cm}\\
b_{2}|_{b_0=0} =&\!\!\!\dfrac{\Gamma(1/4)\sqrt{\pi}}{4\Gamma(11/4)(-2h^2_{20}-2h_{40})^{5/4}}[\delta_{21}(c_{20}-h_{11}c_{11})+\delta_{22}c_{20}],\vspace{0.2cm}\\
b_3|_{b_0=0}=&\!\!\!\Gamma(7/4)\sqrt{\pi}\;[\delta_{31}(c_{20}-h_{11}c_{11})+\delta_{32}c_{20}],\vspace{0.2cm}
\end{array}
\end{equation}
where
$$
\begin{array}{rl}
\delta_{21}=&\!\!\!h^2_{11}-2h_{20},\;\delta_{22}=h^4_{11}-3h^2_{11}h_{20}-2h^2_{20}-4h_{40}\,,\vspace{0.2cm}\\
\delta_{31}=&\!\!\!3h^4_{11}-12h^2_{11}h_{20}+4h^2_{20}-8h_{40}\,,\vspace{0.2cm}\\
\delta_{32}=&\!\!\!3h^6_{11}-15h^4_{11}h_{20}+8h^2_{11}(h^2_{20}-2h_{40})+12h^3_{20}+24h_{20}h_{40}\,.\vspace{0.2cm}
\end{array}
$$

Let
\begin{equation*}
\widetilde{b}_j=b_j|_{b_0=0},\ \ j=1,2,3.
\end{equation*}
Noticing \eqref{b1111}, we can obtain $c_{20}=h_{11}c_{11}-2c_{02}(h^2_{11}-h_{20})$ by solving $\widetilde{b}_1=0$.
Substituting $c_{20}=h_{11}c_{11}-2c_{02}(h^2_{11}-h_{20})$ into $\widetilde b_2$ and $\widetilde b_3$
results in
$$\begin{array}{rl}
\widetilde{b}_2=&\!\!\! \dfrac{\Gamma(1/4)\sqrt{\pi}}{\Gamma(11/4)(-2h^2_{20}-2h_{40})^{1/4}}c_{02},\vspace{0.2cm}\\
\widetilde{b}_3
=&\!\!\!-16\Gamma(7/4)\sqrt{\pi}\;(h^2_{11}-2h_{20})(h^2_{20}+h_{40})c_{02}.\vspace{0.2cm}
\end{array}
$$

For $(h_{20},h_{40})$ near $(h^*_{20},h^*_{40})$, we have $-2h^2_{20}-2h_{40}>0$. Then $\widetilde{b}_2=0$ if and only if
$c_{02}=0$,
which yields $\widetilde{b}_3=0$. Hence, it is clear that the conclusion follows from Theorem 3.

The proof is complete.
\vspace{0.3cm}

Therefore, system \eqref{sns} has $2$ limit cycles near the nilpotent center $(0,1)$. For the symmetry, system \eqref{sns} also has 2 limit cycles near the nilpotent center $(0,-1)$. Namely, there are altogether 4 limit cycles of system \eqref{sns} bifurcating from both symmetric nilpotent centers $(0,\pm1)$.
\vspace{0.3cm}


{\bf Appendix}

Basing on the formulae given in section 2 we have written the Mathematica
code for computing $\{b_j\}$, presented below. It contains
several subroutines (as shown in the code) for computing $\, e_j$,
$q_j(x)$, $H^*_j(x)$, $a_j(x)$, $\bar{q}_j(x)$, $\psi(x)$,
$\hat{q}_j$, $r_{ij}$, $\beta_{ij}$ and $ b_j$.

\footnotesize
\renewcommand{\baselinestretch}{0.8}
\begin{verbatim}
##########  compute the  e_j=e[j]  coefficients  ##########
H=h[2,0]x^2+h[0,2]y^2+Sum[Sum[h[i-j,j]x^{i-j}y^j,{j,0,i}],{i,3,n+1}];
div=Sum[Sum[c[i-j,j]x^{i-j}y^j,{j,0,i}],{i,0,n-1}];
phi=Sum[e[j]x^j,{j,2,m}];
H1=D[H,y];H1=H1/.{y->phi};
Do[H2[i]=Coefficient[H1,x,i],{i,2,m}];
Do[t[i]=Solve[H2[i]==0,e[i]],{i,2,m}];
?t
##########  compute  q_j(x)=q[j]  ##########
Do[q[j]=(1/j!)D[div,{y,j-1}],{j,1,n}];
Do[q[j]=q[j]/.{y->phi},{j,1,n}]
##########  compute H*_j(x)=Hs[j]  ##########
Do[Hs[j]=(1/(j+1)!)D[div,{y,j+1}],{j,1,n}];
Do[Hs[j]=Hs[j]/.{y->phi},{j,1,n}]
##########  compute  a_j(x)=a[j]  ##########
V1=Sum[a[j]w^j,{j,1,m}]
Ht=Sum[Hs[j]V^{j+1},{j,1,n}]
F=Ht-w^2;
Do[W[j]=Coefficient[F,w,j],{j,1,m}]
Do[t1[j]=Solve[W[j+1]==0,a[j]],{j,1,m}]
?t1
##########  compute  bar_q_j(x)=qb[j]  ##########
Q=Sum[q[j]V1^j,{j,1,n}]
Do[qb[j]=2Coefficient[Q,w,2j+1],{j,0,m}]
###  compute  psi(x)=psi  ###
H0=H/.{y->phi}
psi=[H0]^{1/(2p)};
##########  compute  tilde_q_j=qt[j]  ##########
Psi=Sum[PS[j]u^j,{j,1,m}]
H00=H0/.{x->Psi}
G=H00^{2p}-u
Do[g[j]=Coefficient[G,u,j],{j,0,m}]
Do[t2[j]=Solve[g[j]==0,PS[j]],{j,1,m}]
?t2
Do[qt[j]=(qb[j]/.{x->Psi})D[Psi,u],{j,0,m}]
##########  compute  r_ij  coefficients  ##########
Do[r[i,j]=2Coefficient[qt[j],2i],{i,0,m},{j,0,m}]
##########  compute  beta_ij=bt[i,j]  coefficients  ##########
Do[bt[i,j]=Integrate[u^{(1-p+2i)/p}(1-u^2)^j Sqrt[1-u^2],{u,0,1}],{i,0,m},{j,0,m}]
##########  compute  b_j  coefficients  ##########
Do[b[j]=Sum[r[j-p*i,i]bt[j-p*i,i],{i,0,IntegerPart[j/p]}],{j,0,m}]
\end{verbatim}
\normalsize \baselineskip16pt

Executing the above program yields expressions in the
original coefficients $\, h_{ij}$ and $\, c_{ij}$. Below we list the final coefficients $b_j$ (the
intermediate expressions such as $e_j$, $a_j$, etc. are omitted) for $p=2$, $\omega=1$ and $n=3$. We found
$$
b_0=\dfrac{2\pi}{(4 h_{40}-2h^2_{21})^{1/4}}\,c_{00}.
$$
If $\, b_0 =0$, i.e., $\, c_{00} = 0$, then
$$
\begin{array}{rl}
b_1=&\!\!\!\dfrac{16}{3\ (4 h_{40}-2h^2_{21})^{7/4}}\Big[(4 h_{40}-2h^2_{21})c_{20}-(h_{12}h^2_{21}-3h_{21}h_{31}+4h_{12}h_{40})c_{10}\vspace{0.14cm}\\
&\!\!\!-h_{21}(4 h_{40}-2h^2_{21})c_{01}\Big],  \vspace{0.14cm}\\
b_2 = &\!\!\!\dfrac{\pi}{(4 h_{40}-2h^2_{21})^{17/4}} \times
\Big\{2(h_{21}^2 - 2h_{40})\Big[-12h_{21}^4h_{22} + 30h_{12}h_{21}h_{31}(h_{21}^2 + 4h_{40})\vspace{0.15cm}\\
&\!\!\!\;\; + h_{21}^2(-35h_{31}^2+ 8h_{22}h_{40})+4h_{40}(-5h_{31}^2 + 8h_{22}h_{40}) + 3h_{12}^2(h_{21}^4 - 24h_{21}^2h_{40} - 16h_{40}^2)\vspace{0.15cm}\\
&\!\!\!\;\;-4h_{03}h_{21}(h_{21}^4 - 14h_{21}^2h_{40} + 24h_{40}^2)\Big]c_{20}\vspace{0.15cm}\\
&\!\!\!-32 (h_{21}^2 - 2 h_{40})^3 h_{40}c_{02}\vspace{0.15cm}\\
&\!\!\!-8 (h_{21}^2 -2 h_{40})^2 \Big[h_{12} h_{21} (h_{21}^2 - 12 h_{40}) +
   h_{31} (3 h_{21}^2 + 4 h_{40})\Big]c_{11}\vspace{0.15cm}\\
&\!\!\!+2 (h_{21}^2 -2 h_{40})\Big[ 4 h_{03} (h_{21}^2 - 2 h_{40}) (h_{21}^4 - 24 h_{40}^2) +
   6 h_{12} h_{31} (h_{21}^4 - 24 h_{21}^2 h_{40} - 16 h_{40}^2) \vspace{0.15cm}\\
&\!\!\!\;\; + 5 h_{12}^2 h_{21} (h_{21}^4 - 8 h_{21}^2 h_{40} + 48 h_{40}^2) +
   h_{21} (-4 h_{21}^4 h_{22} +      12 h_{40} (5 h_{31}^2 - 8 h_{22} h_{40}) \vspace{0.15cm}\\
&\!\!\!\;\;+      h_{21}^2 (15 h_{31}^2 + 56 h_{22} h_{40}))\Big]c_{01}\vspace{0.15cm}\\
&\!\!\!-8 h_{13} h_{21}(h_{21}^2 - 12 h_{40}) \Big[(h_{21}^2 - 2 h_{40})^2 -
 15 h_{12}^2 h_{21} h_{31} (h_{21}^4 - 40 h_{21}^2 h_{40} - 80 h_{40}^2) \vspace{0.15cm}\\
&\!\!\!\;\;-5 h_{12}^3 (h_{21}^6 - 12 h_{21}^4 h_{40} +
    144 h_{21}^2 h_{40}^2 + 64 h_{40}^3) +
 3 h_{31} (4 h_{03} (h_{21}^6 - 26 h_{21}^4 h_{40} + 32 h_{21}^2 h_{40}^2 \vspace{0.15cm}\\
&\!\!\!\;\;+
       32 h_{40}^3)+    5 h_{21} (4 h_{21}^4 h_{22} +
       4 h_{40} (3 h_{31}^2 - 8 h_{22} h_{40}) +
       h_{21}^2 (7 h_{31}^2 + 8 h_{22} h_{40}))) \vspace{0.15cm}\\
&\!\!\!+ h_{12} (20 h_{03} h_{21} (h_{21}^2 - 2 h_{40}) (h_{21}^4 - 8 h_{21}^2 h_{40} +
       48 h_{40}^2) +
    3 (4 h_{21}^6 h_{22} +
       16 h_{40}^2 (-5 h_{31}^2 \vspace{0.15cm}\\
&\!\!\! + 8 h_{22} h_{40})+
       8 h_{21}^2 h_{40} (-35 h_{31}^2 + 16 h_{22} h_{40}) -
       h_{21}^4 (35 h_{31}^2 + 104 h_{22} h_{40}))\Big]c_{10}\Big\},\vspace{0.15cm}\\
b_3=&\!\!\!\dfrac{16}{15(4 h_{40}-2h^2_{21})^{27/4}} \times(\cdots)\ \ .
\end{array}
$$

\vspace*{1cm}

\noindent {\bf\Large References}

\begin{description}

\item[] Bautin, N. N. [1952]
``On the number of limit cycles which appear with the variation of
coefficients from an equilibrium position of focus or center type,"
 {\em Mat. Sb (N.S.)} {\bf 30},  181-196.

\vspace{-0.05in}

\item[]  Christopher, C, J. \& Lloyd,  N. G. [1996]
 ``Small-amplitude limit cycles in polynomial Li$\acute{e}$nard
 systems,"
 {\em Nonlinear Differential Equations
Appl.}  {\bf 3}, 183-190.

\vspace{-0.05in}

\item[] Christopher, C. J. \& Lynch, S. [1999]
  ``Small-amplitude limit cycle
bifurcations for Li$\acute{e}$nard systems with quadratic or cubic
damping or restoring forces." {\em Nonlinearity} {\bf 12},
1099-1112.

\vspace{-0.05in}

\item[] Gasull, A. \& Torregrosa, J.  [1999]
 ``Small-amplitude limit cycles
in Li$\acute{e}$nard systems via multiplicity,"  {\em J.
Differential Equations} {\bf 159},  186--211.

\vspace{-0.05in}

\item[]Gasull, A. \& Torregrosa, J. [2001]
``A new approach to the computation of the Lyapunov constants. The
geometry of differential equations and dynamical systems,"  {\em
Comput. Appl. Math} {\bf 20}, 149-177.

\vspace{-0.05in}

\item[] Han, M. [1999]
 ``Liapunov constants and Hopf cyclicity of Li$\acute{e}$nard
systems,"  {\em Ann. Diff. Eqns} {\bf 15}, 113-126.

\vspace{-0.05in}

\item[] Han, M.  [2000] ``On Hopf cyclicity of planar systems,"" {\em J. Math.
Anal. Appl.}  {\bf 245} 404-422.

\vspace{-0.05in}

\item[] Han, M. [2006]
 ``Bifurcation Theory of Limit Cycles
of Planar Systems," {\em Handbook of Differential Equations,
Ordinary Differential Equations} {\bf  vol. 3} chapter 4  Edited by
A Canada, P  Drabek and A  Fonda,  Elsevier.

\item[] Han, M. [2012]
{\em Bifurcation Theory of Limit Cycles, Monograph} (Science Press, Beijing, China).

\vspace{-0.05in}

\item[] Han, M.  Chen,  G.  \&  Sun,  C. [2007]
 ``On the number of limit cycles in  near-Hamiltonian polynomial systems,"  {\em Int. J.
Bifurcation and Chaos} {\bf 17}, 2033-2047.

\vspace{-0.05in}

\item[] Han, M.  Lin, Y.  \& Yu, P. [2004]
 ``A study on the existence of limit cycles of a planar
system with third-degree polynomials," {\em Int. J. Bifurcation and
Chaos} {\bf 14}, 41-60.

\vspace{-0.05in}

\item[] Han, M.  Hong, Z.  \& Yang, J. [2009]
 ``Limit cycle bifurcations by perturbing a cuspidal loop in a Hamiltonian system," {\em J. Differential Equations} {\bf 246}, 129-163.

\vspace{-0.05in}

\item[] Han, M.  Jiang, J.  \& Zhu, H. [2008]
 ``Limit cycle bifurcations in near-Hamiltonian system by perturbing a nilpotent center," {\em Int. J. Bifurcation and
Chaos} {\bf 18}, 3013-3027.

\vspace{-0.05in}

\item[] Han, M.  Yang, J.  \& Yu, P. [2009]
 ``Hopf bifurcations for near-Hamiltonian systems," {\em Int. J. Bifurcation and
Chaos} {\bf 19}, 4117-4130.

\vspace{-0.05in}

\item[] Han, M.  Shu, C. Yang, J. \& Chian, A.C.L. [2010]
 ``Polynomial Hamiltonian systems with a nilpotent critical point," {\em Advances in Space Research} {\bf 46}, 521-525.

\vspace{-0.05in}

\item[] Hou, Y. \&  Han, M.  [2006]
``Melnikov functions for planar near-Hamiltonian systems and Hopf
bifurcations," {\em  J. Shanghai Normal University(Natural
Sciences)} {\bf 35},  1-10.

\vspace{-0.05in}

\item[] James, E.  M. \& Lloyd, N. G.  [1991]
``A cubic system with eight small-amplitude limit cycles," {\em IMA
J. Appl. Math} {\bf 47},  163-171.

\vspace{-0.05in}

\item[] Li, J. [2003]
``Hilbert's 16th problem and bifurcations of planar polynomial vector
fields," {\em Int. J. Bifurcation and Chaos} {\bf 13}, 47-106.

\vspace{-0.05in}

\item[]  Llibre, J. [2004]  ``Integrability of Polynomial
Differential Systems. Handbook of Differential Equations," {\em
Ordinary Differential Equations}  {\bf vol. 1}  chapter 5, Edited by
A Canada, P Drabek and A Fonda, Elsevier.

\vspace{-0.05in}

\item[] Li, W.,  Llibre, J. \& Zhang,  X.  [2004]
``Melnikov functions for period annulus, nondegenerate centers,
heteroclinic and homoclinic cycles," {\em Pacific J. of Mathematics}
{\bf 213}, 49-77.

\vspace{-0.05in}

\item[] Schlomiuk, D.  [1993]
``Algebraic and geometric aspects of the theory of polynomial vector
fields, in Bifurcations and Periodic Orbits of Vector Fields. ed.
Schlomiuk D.," {\em NATO ASI Series C} {\bf  408} (Kluwer Academic
London), 429-467.

\vspace{-0.05in}

\item[] Romanovski, V.G. \&  Shafer D.S. [2009]  {\em The Center and Cyclicity Problems: A computational algebra approach} (Birkhauser Boston, Inc., Bostom, MA).

\vspace{-0.05in}

\item[] Yu, P. \&  Han, M. [2004]  ``Twelve limit cycles in a cubic order planar
system with $Z_2$-symmetry," {\em  Communications on Pure and
Applied Analysis} {\bf 3},  515-526.

\end{description}

\end{document}